\documentclass[12pt]{article}
\usepackage{amsmath,amssymb,amscd,latexsym}

\newtheorem{Theorem}{Theorem}[section]

\newtheorem{Proposition}[Theorem]{Proposition}
\newtheorem{Lemma}[Theorem]{Lemma}
\newtheorem{Corollary}[Theorem]{Corollary}

\newtheorem{Remark}[Theorem]{Remark}

\begin{document}

\title{Landis-Oleinik Conjecture in the Exterior Domain}

\author{Jie Wu \,and\, Liqun Zhang\thanks{\small The research is partially
supported by the Chinese NSF under grant 11471320,
 the innovation program at CAS and National Basic Research Program of China under grant 2011CB808002.}}

\date{June, 2015}

\maketitle

\begin{abstract}
In 1974, Landis and Oleinik conjectured that if a bounded solution
of a parabolic equation decays fast at a time, then the solution
must vanish identically before that time, provided the coefficients
of the equation satisfy appropriate conditions at infinity. We prove
this conjecture under some reasonable assumptions on the
coefficients which improved the earlier results.
\end{abstract}



\section{Introduction}
The behavior of solutions of heat equations arose many interests in
last few decades.
In 1974, Landis and Oleinik \cite{LO} proposed the following conjecture:\\

\emph{ If $u(x,t)$ is a bounded solution of a uniformly parabolic equation
$$\sum_{i,j}\partial_i(a^{ij}(x)\partial_ju)-\partial_tu+
b(x)\cdot\nabla u+c(x)u=0~~~~in~~\mathbb{R}^n\times[0,T],
$$
and the condition
\begin{equation}\label{fastdecay}
u(x,T)\leq Ne^{-|x|^{2+\varepsilon}},~~~x\in\mathbb{R}^n,
\end{equation}
holds for some positive constants $N$ and $\varepsilon$,
then $u(x,t)\equiv0$ in $\mathbb{R}^n\times[0,T]$,
provided that the coefficients of the equation satisfy appropriate conditions at infinity.}

The original conjecture only assumes that the coefficients are
time-independent and does not mention the precise conditions,
however, the Lipschitz continuous assumption with some decay at
infinity on $a^{ij}(x)$ seems reasonable and we may also consider
the \emph{space-time} dependent case.

Landis-Oleinik conjecture is closely related to many important
problems. In particular, if $u(x,T)=0$, the conjecture is reduced to
the backward uniqueness problem for parabolic equations. The
backward uniqueness problem has a natural background in the control
theory for PDEs, and it also appeared in the regularity theory of
parabolic equations, such as the Navier-Stokes equations
\cite{ESS1}, semi-linear heat equations \cite{Miz}, heat flow of
harmonic maps \cite{Wa}.

This conjecture has a elliptic version, where probably the problem
originated, the Landis conjecture, namely, if a solution of an
elliptic equation decays faster than a given rate at infinity, then
it is identically zero. The complex case of Landis conjecture is
solved by Meshkov \cite{Me}, and a quantitative result is proved by
Bourgain and Kenig \cite{BK}, while the real case remains open.

Now we denote the backward parabolic operator
$$P=\partial_t+\sum_{i,j}\partial_i(a^{ij}(x,t)\partial_j)=\partial_t+\nabla\cdot(\mathbf{A}\nabla),$$
where $\mathbf{A}(x,t)=(a^{ij}(x,t))^n_{i,j=1}$ is a real symmetric matrix
such that for some $\Lambda\geq\lambda>0$,
\begin{equation}\label{a1}
\lambda|\xi|^2\leq \sum_{i,j}a^{ij}(x,t)\xi_i\xi_j\leq
\Lambda|\xi|^2, ~\forall \xi\in\mathbb{R}^n.
\end{equation}
In the following we always assume that the lower coefficients $b$
and $c$ are space-time dependent and bounded,
and rather than (\ref{fastdecay}), we assume
a weaker condition
\begin{equation}\label{fastdecay1}
|u(x,0)|\leq C_ke^{-k|x|^2},~~~\forall~~ k>0.
\end{equation}

There are some earlier results about Landis-Oleinik conjecture. In
the constant coefficients case, i.e., the heat equation, this
conjecture was solved by Escauriaza, Kenig, Ponce and
Vega \cite{EKPV}. They introduced some interesting Carleman
estimates and proved both qualitative and quantitative results
for the heat equation with bounded space-time dependent coefficients
in whole space and half space.

For the general case, the first result is obtained by Nguyen \cite{Tu} where  both
qualitative and quantitative results are proved for the conjecture
in whole space and half space under the following assumptions
\begin{equation}\label{Nguyen1}
|\nabla_xa^{ij}(x,t)|+|\partial_ta^{ij}(x,t)|\leq M,
\end{equation}
\begin{equation}\label{Nguyen2}
|\nabla_xa^{ij}(x,t)|\leq M\langle x\rangle^{-1-\epsilon},
\end{equation}
\begin{equation}\label{Nguyen3}
|a^{ij}(x,t)-a^{ij}(x,s)|\leq M\langle x\rangle^{-1}|t-s|^{1/2},
\end{equation}
where $\langle x\rangle=\sqrt{1+|x|^2}$ and $\varepsilon>0$. We
remark that condition (\ref{Nguyen1}), the Lipschitz regularity
assumption is reasonable, as shown in \cite{Pl,Mil}, and some decay
assumptions seems necessary. However, condition (\ref{Nguyen2}) is
not scaling invariant and we wonder if condition (\ref{Nguyen3}) is necessary.

Another related result is
the backward uniqueness result for general parabolic equations in
half space proved by the authors \cite{WZ} under condition
(\ref{Nguyen1}) and the decay at infinity condition:
\begin{equation}\label{opcond}
|\nabla_xa^{ij}(x,t)|\leq \frac{E}{|x|},~~~where~~E<E_0(n,\Lambda,\lambda).
\end{equation}
Note that condition (\ref{opcond}) is scaling invariant. In \cite{WZ},
the authors also constructed examples
to show that both condition (\ref{Nguyen1}) and (\ref{opcond}) are almost optimal.

All these results suggest that if Landis-Oleinik conjecture is true,
then certain regularity and decay at infinity assumptions on the
coefficients should be required, and assumptions (\ref{Nguyen1}) and
(\ref{opcond}) seem to be optimal.

Now in the exterior domain, under assumptions (\ref{Nguyen1}) and
(\ref{opcond}), we shall prove the Landis-Oleinik conjecture. Our
main result is the following.
\begin{Theorem}\label{mainthm}
Suppose $\{a^{ij}\}$ satisfy (\ref{a1}), and for some constants $E,M,N>0$,
\begin{equation}\label{a2}
|\nabla_xa^{ij}(x,t)|+|\partial_ta^{ij}(x,t)|\leq M, ~~~~
\forall ~(x,t)\in \mathbb{R}^n\backslash B_1\times[0,1],
\end{equation}
and
\begin{equation}\label{a3}
|\nabla_xa^{ij}(x,t)|\leq \frac{E}{|x|},  ~~~~~~~~\forall ~(x,t)\in
\mathbb{R}^n\backslash B_1\times[0,1].
\end{equation}
Assume that $u$ satisfies
\begin{eqnarray}\label{equat}
\left\{
  \begin{array}{ll}
    |Pu|\leq N(|u|+|\nabla u|)~~~~&in~~ \mathbb{R}^n\backslash B_1\times[0,1], \\
    |u(x,t)|\leq Ne^{N|x|^2}~~~~&in~~ \mathbb{R}^n\backslash B_1\times[0,1],\\
    |u(x,0)|\leq C_ke^{-k|x|^2},~~~\forall k>0~~~~&in~~\mathbb{R}^n\backslash B_1.
  \end{array}
\right.
\end{eqnarray}
Then there exists a constant $E_0=E_0(n,\Lambda,\lambda)$, such that when $E<E_0$, we have\\
%%%  1. if $u(x,0)\not\equiv0$, there exists $N>0$ such that if $|x|>N,$
%%%  $$\int_{B(x,1)}|u(y,0)|^2dy\geq e^{-N|x|^2log|x|}~~~and~~~\int_{B(x,\frac{|x|}{2})}|u(y,0)|^2dy\geq e^{-N|x|^2}.$$
%%%  2. if $|u(x,0)|\leq C_ke^{-k|x|^2}$ for all $k>0$, $x\in \mathbb{R}^n\backslash B_R$,
$u(x,t)\equiv 0$ in $\mathbb{R}^n\backslash B_1\times[0,1]$.
\end{Theorem}

By the unique continuation(see \cite{EF,F}) result, we immediately have the following corollary.
\begin{Corollary}
Theorem \ref{mainthm} is still valid if we replace $\mathbb{R}^n\backslash B_1$ by $\mathbb{R}^n$.
\end{Corollary}

Theorem \ref{mainthm} can be obtained immediately by the following upper bound and lower bound estimates.

\begin{Proposition}[Upper Bound]\label{ubound}
Suppose $\{a^{ij}\}$ and $u$ are the same as above.
Then there exists a constant $E_0=E_0(n,\Lambda,\lambda)$, such that when $E<E_0$, we have
$$|u(x,t)|+|\nabla u(x,t)|\leq e^{-k|x|^2},~~~~~~~\forall k>0,$$
when $|x|\geq R_1(n,\Lambda,\lambda,M,E,N,k)$ and $0\leq t\leq
T_1(\Lambda,N)$.
\end{Proposition}

\begin{Proposition}[Lower Bound]\label{lbound}
Suppose $\{a^{ij}\}$ are the same as above, $u$ satisfies the first
two conditions of (\ref{equat}), and $u(x,0)\neq0$. Then there
exists a positive constant $E_0=E_0(n, \Lambda, \lambda)$, such that
when $E< E_0$, there exists a constant
$C_\star=C_\ast(n,\Lambda,\lambda,M,E,N)$, such that the following
estimate
\begin{equation}\label{lowerbound}
\frac{1}{T}\int_{T/8}^{7T/8}\int_{R-1\leq|x|\leq R}(u^2+|\nabla u|^2)dxdt\geq e^{-C_\star\frac{R^2}{T}}.
\end{equation}
holds when $$R\geq
R_2(n,\Lambda,\lambda,M,E,N,||u(\cdot,0)||_{L^2(B(10e_1,\frac{1}{2}))})$$
and $$0<T\leq
T_2(n,\Lambda,\lambda,M,N,||u(\cdot,0)||_{L^2(B(10e_1,\frac{1}{2}))}),$$
where $e_1=(1,0,\ldots,0)$.
\end{Proposition}

Combining these two estimates together, we must have $u(x,0)=0$,
then by the backward uniqueness(see \cite{WZ}) result, we have
$u(x,t)\equiv0$. Thus we proved Theorem \ref{mainthm}.

\begin{Remark}
This lower bound of the integration form is optimal, which can be seen from
the solution of the backward heat equation $\partial_t\Gamma+\triangle\Gamma=0$ that
$$\Gamma(x,t)=(T-t)^{-n/2}e^{-\frac{|x|^2}{4(T-t)}}.$$
\end{Remark}

The upper bound can be obtained by the following Carleman inequality.
\begin{Proposition}\label{Prop-C1}
Suppose $\{a^{ij}\}$ are the same as above. Let \\
$$Q=\mathbb{R}^n\backslash B_1\times[0,1],~~~~f(t)=(t+1)^{-\beta}-2^{-\beta}.$$
There exists a positive constant $E_0=E_0(n, \Lambda, \lambda)$, such that when $E< E_0$, for any
function $v\in C^\infty_0(Q)$ and any $\gamma>0$, we have
\begin{equation}\label{C1}
\begin{split}
&\int_{Q}e^{2\gamma f|x|^{3/2}-\frac{b|x|^2+\beta}{t+1}}(|v|^2+|\nabla v|^2)dxdt
\leq \int_{Q}e^{2\gamma f|x|^{3/2}-\frac{b|x|^2+\beta}{t+1}}|Pv|^2dxdt\\
&~~~~~~~~~~~~~~+c\int_{\mathbb{R}^n}|x|^2 e^{-\frac{b|x|^2}{2}}(|v(x,1)|^2+|\nabla v(x,1)|^2)dx\\
&~~~~~~~~~~~~~~+c(1+\gamma)^2\int_{\mathbb{R}^n}|x|^2e^{2\gamma|x|^{3/2}-b|x|^2}|v(x,0)|^2dx,
\end{split}
\end{equation}
where $b=\frac{1}{16\Lambda}$, $\beta=\beta(n,\Lambda,\lambda,M,E)\geq1$, and $c$ is an absolute constant.
\end{Proposition}

The lower bound can be proved  mainly by the following Carleman inequality.\\
First, let $\psi(t)$ be a cut-off function satisfying
$$
\psi(t)=\left\{
                 \begin{array}{ll}
                   0, & \hbox{if $t\in[0,\frac{1}{4}]\cup[\frac{3}{4},1];$} \\
                   2, & \hbox{if $t\in[\frac{1}{3},\frac{2}{3}]$.}
                 \end{array}
               \right.
$$

\begin{Proposition}\label{Prop-C2}
Suppose $\{a^{ij}\}$ are the same as above. Let
$$Q_R=\{(x,t)|\quad 1<|x|<R, t\in(\frac{1}{8},\frac{7}{8})\},$$
$$\Psi=\gamma(1-t)R^{2/3}|x|^{4/3}+\psi(t)R^2.$$
Then there exists a positive constant $E_0=E_0(n, \Lambda, \lambda)$,
such that when $E< E_0$, for any
function $v\in C^\infty_0(Q_R)$ and any $\gamma\geq\gamma_0(n,\Lambda,\lambda,M,E)$, we have
\begin{equation}\label{C2}
c\lambda^2\int_{Q_R}e^{2\Psi}(\gamma^3R^2|v|^2+\gamma|\nabla v|^2)dxdt\leq \int_{Q_R}e^{2\Psi}|Pv|^2dxdt.
\end{equation}
\end{Proposition}

The paper organized as follows. We first use the two Carleman
inequalities to prove the upper and lower bound, then we prove the two Carleman inequalities.

\section{Proof of Upper Bound and Lower Bound}

In this section, we prove the upper bound and lower bound by
assuming Proposition \ref{Prop-C1} and Proposition \ref{Prop-C2}
first, and we postpone the proof of the two Carleman inequalities to
the next section.

\subsection{Upper Bound}

\emph{Proof of Proposition \ref{ubound}}. We use Carleman inequalities (\ref{C1})
to prove the upper bound for the solution.\\

\textbf{Step 1.} By the regularity theory for solutions of parabolic equations, we have
\begin{equation}\label{pt11}
|u(x,t)|+|\nabla u(x,t)|\leq C(n,\Lambda,\lambda,M,N)e^{2N|x|^2}
\end{equation}
for $(x,t)\in \mathbb{R}^n\backslash B_2\times[0,\frac{1}{2}]$.
Let
\begin{equation}\label{T_1variable}
\tau=\min\{\frac{1}{2}, \frac{1}{2N},\frac{b}{8}\},
\end{equation}
where $b$ is the one in Proposition \ref{Prop-C1}. Define
$$\tilde{u}(x,t)=u(\tau x,\tau^2t),$$
and
$$\tilde{a}^{ij}(x,t)={a}^{ij}(\tau x,\tau^2t)$$
for $(x,t)\in\mathbb{R}^n\backslash B_{\frac{2}{\tau}}\times[0,1]$. Then it is easy to see
$$|\nabla_x\tilde{a}^{ij}|+|\partial_t\tilde{a}^{ij}|
\leq \tau M\leq M,~~~|\nabla_x\tilde{a}^{ij}|\leq\frac{E}{|x|}.$$
We denote
$$\tilde{P}\tilde{u}=\partial_t\tilde{u}+\sum_{ij}\partial_i(\tilde{a}^{ij}\partial_j\tilde{u}),$$
then by (\ref{equat}) we have
\begin{equation}\label{pt2}
|\tilde{P}\tilde{u}|\leq\tau N(|\tilde{u}|+|\nabla \tilde{u}|)
\leq\frac{1}{2}(|\tilde{u}|+|\nabla \tilde{u}|).
\end{equation}
By (\ref{pt11}) and (\ref{T_1variable}), we have
\begin{equation}\label{pt12}
|\tilde{u}(x,t)|+|\nabla \tilde{u}(x,t)|\leq C(n,\Lambda,\lambda,M,N)e^{2N\tau^2|x|^2}
\leq C(n,\Lambda,\lambda,M,N)e^{\frac{b}{8}|x|^2}.
\end{equation}
We keep in mind that
\begin{equation}\label{initialtime}
|u(x,0)|\leq C_ke^{-k|x|^2},~~~\forall k>0,
\end{equation}
and we always take $k$ large enough.\\

\textbf{Step 2.}
In order to apply Carleman inequality (\ref{C1}), we define a cut-off function $\theta$ satisfying
$$
\theta(|x|)=\left\{
                 \begin{array}{ll}
                   0, & \hbox{if $|x|<R~~or~~|x|>k^2R+1$}; \\
                   1, & \hbox{if $R+1\leq|x|\leq k^2R$},
                 \end{array}
               \right.
$$
where $R>\frac{2}{\tau}$.\\
Let $v=\theta\tilde{u}$, then by (\ref{pt2}) we have
\begin{equation}\label{estv}
\begin{split}
|\tilde{P}v|=&|\theta\tilde{P}\tilde{u}+\tilde{u}\tilde{P}\theta+
2\tilde{a}^{ij}\partial_i\theta\partial_j \tilde{u}|\\
\leq&\frac{1}{2}\theta(|\tilde{u}|+|\nabla \tilde{u}|)+
C(n,\Lambda,M)\chi(|\tilde{u}|+|\nabla \tilde{u}|)(|\nabla\theta|+|\nabla^2\theta|)\\
\leq&\frac{1}{2}(|v|+|\nabla v|)+C(n,\Lambda,M)\chi_{\Omega}(|\tilde{u}|+|\nabla \tilde{u}|),
\end{split}
\end{equation}
where $\chi$ is the characteristic function and
\begin{eqnarray*}
\begin{split}
\Omega=&\{|0<\theta<1, t\in[0,1]\}\\
=&\{R<|x|<R+1, t\in[0,1]\}\bigcup\{k^2R<|x|<k^2R+1, t\in[0,1]\}.
\end{split}
\end{eqnarray*}

\textbf{Step 3.} We apply Carleman inequality (\ref{C1}) for $v$, then
\begin{eqnarray*}
\begin{split}
J\equiv&\int_{Q}e^{2\gamma f|x|^{3/2}-\frac{b|x|^2+\beta}{t+1}}(|v|^2+|\nabla v|^2)dxdt\\
\leq &\int_{Q}e^{2\gamma f|x|^{3/2}-\frac{b|x|^2+\beta}{t+1}}|\tilde{P}v|^2dxdt\\
&~+c\int_{\mathbb{R}^n}|x|^2 e^{-\frac{b|x|^2}{2}}(|v(x,1)|^2+|\nabla v(x,1)|^2)dx\\
&~+c(1+\gamma)^2\int_{\mathbb{R}^n}|x|^2e^{2\gamma|x|^{3/2}-b|x|^2}|v(x,0)|^2dx.\\
\end{split}
\end{eqnarray*}
By (\ref{estv}) we have
\begin{eqnarray*}
\begin{split}
J\leq &\frac{3}{4}J+C(n,\Lambda,M)\int_{Q}e^{2\gamma f|x|^{3/2}-
\frac{b|x|^2+\beta}{t+1}}\chi(|\tilde{u}|+|\nabla \tilde{u}|)^2dxdt\\
&~+c\int_{\mathbb{R}^n}|x|^2 e^{-\frac{b|x|^2}{2}}(|v(x,1)|^2+|\nabla v(x,1)|^2)dx\\
&~+c(1+\gamma)^2\int_{\mathbb{R}^n}|x|^2e^{2\gamma|x|^{3/2}-b|x|^2}|v(x,0)|^2dx,\\
\end{split}
\end{eqnarray*}
thus
\begin{eqnarray*}
\begin{split}
J\leq &C(n,\Lambda,M)\int_{\Omega}e^{2\gamma f|x|^{3/2}-\frac{b|x|^2+
\beta}{t+1}}(|\tilde{u}|+|\nabla \tilde{u}|)^2dxdt\\
&~+c\int_{|x|\geq R}|x|^2 e^{-\frac{b|x|^2}{2}}(|\tilde{u}(x,1)|^2+|\nabla \tilde{u}(x,1)|^2)dx\\
&~+c(1+\gamma)^2\int_{|x|\geq R}|x|^2e^{2\gamma|x|^{3/2}-b|x|^2}|\tilde{u}(x,0)|^2dx\\
\equiv&I_1+I_2+I_3.
\end{split}
\end{eqnarray*}

\textbf{Step 4.} Now we estimate both sides of the above inequality.
We estimate $I_2$ and $I_3$ first, then $I_1$, at last $J$.\\

\textbf{Estimate of $I_2$.}

By (\ref{pt12}),
\begin{equation}\label{I2}
I_2\leq C(n,\Lambda,\lambda,M,N)\int_{|x|\geq R}|x|^2e^{-\frac{b|x|^2}{4}}dx\leq C(n,\Lambda,\lambda,M,N)
\end{equation}

\textbf{Estimate of $I_3$.}

Recall (\ref{initialtime}), then
$$|\tilde{u}(x,0)|=|u(\tau x,0)|\leq C(\frac{bk}{\tau^2})e^{-\frac{bk}{\tau^2}|\tau x|^2}
=C(\Lambda,N,k)e^{-bk|x|^2},$$
and thus
$$I_3\leq C(\Lambda,N,k)(1+\gamma)^2\int_{|x|\geq R}|x|^2e^{-bk|x|^2+2\gamma|x|^{3/2}-b|x|^2}dx.$$
Now we choose
\begin{equation}\label{chgamma}
\gamma=\frac{bk}{16}R^{1/2}.
\end{equation}
In the region $\{|x|\geq R\}$,
$$2\gamma|x|^{3/2}=\frac{bk}{8}R^{1/2}|x|^{3/2}\leq\frac{bk}{8}|x|^2,$$
then
\begin{equation}\label{I3}
\begin{split}
I_3\leq& C(\Lambda,N,k)k^2R\int_{|x|\geq R}|x|^2e^{-\frac{bk}{2}|x|^2-b|x|^2}dx\\
\leq& C(\Lambda,N,k)k^2Re^{-\frac{bk}{2}R^2}\int_{|x|\geq R}|x|^2e^{-b|x|^2}dx\\
\leq& C(n,\Lambda,N,k)k^2Re^{-\frac{bk}{2}R^2}\leq 1,
\end{split}
\end{equation}
if $R\geq R_0(n,\Lambda,N,k)$ large enough.

\textbf{Estimate of $I_1$.}

$$I_1\leq C(n,\Lambda,M)\int_{\Omega}e^{2\gamma |x|^{3/2}-
\frac{b|x|^2}{2}}(|\tilde{u}|+|\nabla \tilde{u}|)^2dxdt.$$
Use (\ref{pt12}) again, we obtain
\begin{eqnarray*}
\begin{split}
I_1\leq& C(n,\Lambda,\lambda,M,N)\int_{\Omega}e^{2\gamma |x|^{3/2}-\frac{b|x|^2}{4}}dxdt\\
\leq&C(n,\Lambda,\lambda,M,N)(\int_{k^2R\leq|x|\leq k^2R+1}+
\int_{R\leq|x|\leq R+1})e^{2\gamma|x|^{3/2}-\frac{b|x|^2}{4}}dx\\
\equiv & I_{1,1}+I_{1,2}.
\end{split}
\end{eqnarray*}
In the region $\{k^2R\leq|x|\leq k^2R+1\}$,
$$2\gamma|x|^{3/2}=\frac{bk}{8}R^{1/2}|x|^{3/2}\leq\frac{b}{8}|x|^2,$$
then
$$I_{1,1}\leq C(n,\Lambda,\lambda,M,N)\int_{k^2R\leq|x|
\leq k^2R+1}e^{-\frac{b|x|^2}{8}}dx\leq C(n,\Lambda,\lambda,M,N).$$
In $\{R\leq|x|\leq R+1\}$,
$$2\gamma|x|^{3/2}=\frac{bk}{8}R^{1/2}|x|^{3/2}\leq\frac{bk}{8}|x|^2,$$
then
\begin{eqnarray*}
\begin{split}
I_{1,2}\leq &C(n,\Lambda,\lambda,M,N)\int_{R\leq|x|\leq R+1}e^{\frac{bk}{8}|x|^2-\frac{b|x|^2}{4}}dx\\
\leq &C(n,\Lambda,\lambda,M,N)e^{\frac{bk}{8}(R+1)^2}\int_{R\leq|x|\leq R+1}e^{-\frac{b|x|^2}{4}}dx\\
\leq &C(n,\Lambda,\lambda,M,N)e^{\frac{bk}{2}R^2}.
\end{split}
\end{eqnarray*}
Thus we have
\begin{equation}\label{I1}
\begin{split}
I_1\leq C(n,\Lambda,\lambda,M,N)e^{\frac{bk}{2}R^2}.
\end{split}
\end{equation}
Combining (\ref{I2}), (\ref{I3}) and (\ref{I1}), we have that when
$R\geq R_0(n,\Lambda,N,k)$,
\begin{equation}\label{uboundJ}
J\leq C(n,\Lambda,\lambda,M,N)e^{\frac{bk}{2}R^2}.
\end{equation}

Next we estimate a lower bound for $J$.\\

\textbf{Estimate of $J$.}

If $k\geq 4^{\beta+5}$, then $\{4^{\beta+2}R\leq|x|\leq 4^{\beta+3}R\}\subset\{\theta=1\}$, and thus
$$
J\geq\int_0^{1/2}\int_{4^{\beta+2}R\leq|x|\leq 4^{\beta+3}R}e^{2\gamma f|x|^{3/2}-
\frac{b|x|^2+\beta}{t+1}}(|\tilde{u}|^2+|\nabla \tilde{u}|^2)dxdt.
$$
Notice that when $t\in [0,\frac{1}{2}]$, $f(t)\geq f(\frac{1}{2})\geq2^{-\beta-2}$, then
$$J\geq e^{-\beta}\int_0^{1/2}\int_{4^{\beta+2}R\leq|x|
\leq 4^{\beta+3}R}e^{2^{-\beta-1}\gamma
|x|^{3/2}-b|x|^2}(|\tilde{u}|^2+|\nabla \tilde{u}|^2)dxdt.$$
In the
region $\{4^{\beta+2}R\leq|x|\leq 4^{\beta+3}R\}$,
$$2^{-\beta-1}\gamma |x|^{3/2}=2^{-\beta-5}bkR^{1/2}|x|^{3/2}
\geq2^{-\beta-5}bk(4^{-\beta-3}|x|)^{1/2}|x|^{3/2}=4^{-\beta-4}bk|x|^2,$$
then
$$2^{-\beta-1}\gamma |x|^{3/2}-b|x|^2\geq(4^{-\beta-4}k-1)b|x|^2.$$
Notice that $k\geq 4^{\beta+5}$ , then
$$4^{-\beta-4}k-1\geq4^{-\beta-5}k,$$
and
$$2^{-\beta-1}\gamma |x|^{3/2}-b|x|^2\geq4^{-\beta-5}bk|x|^2
\geq4^{-\beta-5}bk(4^{\beta+2}R)^2=4^{\beta-1}bkR^2\geq bkR^2.$$
Thus
\begin{equation}\label{lboundJ}
\begin{split}
J\geq& e^{-\beta}e^{bkR^2}\int_0^{1/2}\int_{4^{\beta+2}R\leq|x|
\leq 4^{\beta+3}R}(|\tilde{u}|^2+|\nabla \tilde{u}|^2)dxdt\\
\geq& \tau^{-n}e^{-\beta}e^{bkR^2}\int_0^{\tau^2/2}
\int_{\tau4^{\beta+2}R\leq|x|\leq \tau4^{\beta+3}R}(|u|^2+|\nabla u|^2)dxdt.
\end{split}
\end{equation}
Combining (\ref{uboundJ}) and (\ref{lboundJ}) together, we have
$$\int_0^{\tau^2/2}\int_{\tau4^{\beta+2}R\leq|x|
\leq \tau4^{\beta+3}R}(|u|^2+|\nabla u|^2)dxdt
\leq  C(n,\Lambda,\lambda,M,E,N)e^{-\frac{bk}{2}R^2}\leq e^{-\frac{bk}{4}R^2}$$
when $R\geq R_0(n,\Lambda,\lambda,M,E,N,k)$.

We replace $\tau4^{\beta+2}R$ by $R$, and let
$$T_1=\frac{\tau^2}{4}=\frac{1}{16}\min\{1, \frac{1}{N^2},\frac{b^2}{16}\},$$
then we obtain
$$\int_0^{2T_1}\int_{R\leq|x|\leq 4R}(|u|^2+|\nabla u|^2)dxdt\leq e^{-CkR^2}.$$
Finally, by the regularity theory for solutions of parabolic
equations, we obtained our upper bound estimate.

\subsection{Lower Bound}

The lower bound can be proved by the following two lemmas. The first
one is due to Escauriaza, Fern¡äandez and Vessella (see \cite{EFV}), and we copy it here.

\begin{Lemma}\label{doubleproperty}
There is a constant $C=C(n,\Lambda,\lambda,M,N)$ such that the inequalities
\begin{equation}
C\log(C\Theta_\rho)\geq 2~~~~and~~~~C\int_{B_{2\rho}}u^2(x,t)dx\geq\int_{B_{\rho}}u^2(x,0)dx
\end{equation}
hold when $0<t\leq\rho^2/C\log(C\Theta_\rho)$ and $0<\rho\leq1$. Here
$$\Theta_\rho=\frac{\int_0^1\int_{B_4}u^2(x,t)dxdt}{\rho^2\int_{B_\rho}u^2(x,0)dx}.$$
\end{Lemma}

The second one is derived from Carleman inequality (\ref{C2}).
\begin{Lemma}\label{glbound}
Suppose $\{a^{ij}\}$ are the same as above, $u$ satisfies the first two conditions of (\ref{equat}).
Then there exists a positive constant $E_0=E_0(n, \Lambda, \lambda)$, such that when $E< E_0$,
there exists $C_\star=C_\ast(n,\Lambda,\lambda,M,E,N)$, such that the following estimate
\begin{equation}\label{lowerbound}
\begin{split}
&e^{C_\ast\frac{R^{2/3}}{T}}\frac{1}{T}\int_{T/3}^{2T/3}\int_{9\leq|x|\leq11}|u|^2dxdt\\
\leq&1+e^{C_\ast\frac{R^2}{T}}\frac{1}{T}\int_{T/8}^{7T/8}\int_{R-1\leq|x|\leq R}(|u|^2+|\nabla u|^2)dxdt,
\end{split}
\end{equation}
holds when $R\geq R_3(n,N)$ and $0<T\leq1$.
\end{Lemma}

In the following, we prove Lemma \ref{glbound} first, then we use the two lemmas
to prove the lower bound.\\

\emph{Proof of Lemma \ref{glbound}}. We use Carleman inequality
(\ref{C2}) to prove Lemma \ref{glbound}. We again divided the proof
into several steps.

{\bf Step 1}. For any $0<T\leq 1$, we define
$$\tilde{u}(x,t)=u(\sqrt{T}x,Tt),$$
$$\tilde{a}^{ij}(x,t)={a}^{ij}(\sqrt{T}x,Tt),$$
$$\tilde{P}\tilde{u}=\partial_t\tilde{u}+\sum_{ij}\partial_i(\tilde{a}^{ij}\partial_j\tilde{u}),$$
for $(x,t)\in (x,t)\in\mathbb{R}^n\backslash B_{\frac{2}{\sqrt{T}}}\times[0,1]$. Similarly, we have
\begin{equation}\label{pt2a}
|\tilde{u}(x,t)|+|\nabla \tilde{u}(x,t)|\leq C(n,\Lambda,\lambda,M,N)e^{2NT|x|^2},
\end{equation}
and
\begin{equation}\label{pt12a}
|\tilde{P}\tilde{u}|\leq \sqrt{T}N(|\tilde{u}|+|\nabla \tilde{u}|)\leq N(|\tilde{u}|+|\nabla \tilde{u}|).
\end{equation}

\textbf{Step 2.} In order to apply Carleman inequality (\ref{C2}), we choose two
smooth cut-off functions. Let
$$
\eta_1(|x|)=\left\{
                 \begin{array}{ll}
                   0, & \hbox{if $|x|\leq\frac{2}{\sqrt{T}}~~or~~|x|\geq \gamma^{-3/4}R$}; \\
                   1, & \hbox{if $\frac{3}{\sqrt{T}}\leq|x|\leq\gamma^{-3/4}R-\frac{1}{\sqrt{T}}$},
                 \end{array}
               \right.
$$
where $\gamma$ and $R$ are the parameters in Carleman inequality (\ref{C2}), and
\begin{equation}\label{gammar}
\gamma^{-3/4}R\geq\frac{20}{\sqrt{T}}.
\end{equation}
We always take both $\gamma$ and $R$ large enough.
Let
$$
\eta_2(t)=\left\{
            \begin{array}{ll}
              0, & \hbox{if $t\in [0,\frac{1}{8}]\bigcup[\frac{7}{8},1]$} ;\\
              1, & \hbox{if $t\in[\frac{1}{4},\frac{3}{4}]$}.
            \end{array}
          \right.
$$
%%%% Both functions take values in $[0,1]$, and
%%%%\begin{equation}\label{cutoff}
%%%%|\theta_1'|\leq C,~~~|\theta_1''|\leq C,~~~|\theta'_2|\leq C(n)/R,~~and~~|\theta''_2|\leq C(n)/R^2
%%%%\end{equation}
Let $\eta(x,t)=\eta_1(|x|)\eta_2(t)$ and $v=\eta \tilde{u}$.
Then $supp~\eta\subset Q_R$ and so $supp~w\subset Q_R$.\\
By (\ref{pt12a}) we have
\begin{equation}\label{v1a}
\begin{split}
|\tilde{P}v|=&|\eta\tilde{P}\tilde{u}+\tilde{u}\tilde{P}\eta+
2\tilde{a}^{ij}\partial_i\eta\partial_j \tilde{u}|\\
\leq&N\eta(|\tilde{u}|+|\nabla \tilde{u}|)+
C(n,\Lambda,M)(|\tilde{u}|+|\nabla \tilde{u}|)(|\partial_t\eta|+|\nabla\eta|+|\nabla^2\eta|)\\
\leq&N(|v|+|\nabla v|)+C(n,\Lambda,M,N)(|\tilde{u}|+|\nabla \tilde{u}|)\chi_{\{0<\eta<1\}},
\end{split}
\end{equation}

{\bf Step 3.} We apply Calman inequality (\ref{C2}) for $v$, then we get
$$
c\lambda^2\int_{Q_R}e^{2\Psi}(\gamma^3R^2|v|^2+\gamma|\nabla v|^2)dxdt
\leq \int_{Q_R}e^{2\Psi}|\tilde{P}v|^2dxdt.
$$
By (\ref{v1a}), we have
\begin{eqnarray*}
\begin{split}
c\lambda^2\int_{Q_R}e^{2\Psi}(\gamma^3R^2|v|^2+\gamma|\nabla v|^2)dxdt
\leq 4N^2\int_{Q_R}e^{2\Psi}(|v|^2+|\nabla v|^2)dxdt\\
+C\int_{\{0<\eta<1\}}e^{2\Psi}(|\tilde{u}|^2+|\nabla\tilde{u}|^2)dx dt.
\end{split}
\end{eqnarray*}
In the above inequality, if we take $\gamma=\gamma(n,\Lambda,\lambda,M,E,N)$
large enough, then the first term of the right hand side
can be absorbed by the term of the left hand side, thus we obtain
$$\int_{Q_R}e^{2\Psi}(|v|^2+|\nabla v|^2)dxdt\leq C\gamma^{-1}
\int_{\{0<\eta<1\}}e^{2\Psi}(|\tilde{u}|^2+|\nabla\tilde{u}|^2)dxdt.$$
Denote that
$$\Omega_1=\{(x,t)|~\frac{9}{\sqrt{T}}\leq|x|\leq\frac{11}{\sqrt{T}},~t\in[\frac{1}{3},\frac{2}{3}]\},$$
then $\Omega_1\subset\{\eta=1\}$ and thus
\begin{equation}\label{ideaa}
\int_{\Omega_1}e^{2\Psi}(|\tilde{u}|^2+|\nabla \tilde{u}|^2)dxdt
\leq C\gamma^{-1}\int_{\{0<\eta<1\}}e^{2\Psi}(|\tilde{u}|^2+|\nabla\tilde{u}|^2)dxdt.
\end{equation}
We divide the set $\{0<\eta<1\}$ into three parts:
$$\{0<\eta<1\}\subset \Omega_2\cup\Omega_3\cup\Omega_4,$$
where
\begin{equation}\label{supp}
\begin{split}
\Omega_2=&\{(x,t)|~\frac{2}{\sqrt{T}}<|x|<\frac{3}{\sqrt{T}},~t\in[\frac{1}{8},\frac{7}{8}]\},\\
\Omega_3=&\{(x,t)|~\frac{3}{\sqrt{T}}<|x|<\gamma^{-3/4}R-
\frac{1}{\sqrt{T}},~t\in[\frac{1}{8},\frac{1}{4}]\cup[\frac{3}{4},\frac{7}{8}]\},\\
\Omega_4=&\{(x,t)|~\gamma^{-3/4}R-
\frac{1}{\sqrt{T}}<|x|<\gamma^{-3/4}R,~t\in[\frac{1}{8},\frac{7}{8}]\}.\\
\end{split}
\end{equation}
If we denote that
$$J_i=\int_{\Omega_i}e^{2\Psi}(|\tilde{u}|^2+|\nabla\tilde{u}|^2)dxdt,~~~i=1,2,3,4,$$
then we rewrite (\ref{ideaa}) as
\begin{equation}\label{ideab}
J_1\leq C\gamma^{-1}(J_2+J_3+J_4).
\end{equation}

{\bf Step 4.} We estimate them respectively. \\

\textbf{Estimate of $J_1$.}

In $\Omega_1$, $\psi(t)=2$, and
$$\Psi\geq\frac{\gamma}{3} R^{2/3}(\frac{9}{\sqrt{T}})^{4/3}+2R^2\geq 6\gamma (\frac{R}{T})^{2/3}+2R^2,$$
then
\begin{equation}\label{lhsa}
\begin{split}
J_1\geq& \exp\{12\gamma (\frac{R}{T})^{2/3}+4R^2\}\int_{\Omega_1}|\tilde{u}|^2dxdt\\
=&T^{-\frac{n}{2}-1}\exp\{12\gamma (\frac{R}{T})^{2/3}+4R^2\}
\int_{T/3}^{2T/3}\int_{9\leq|x|\leq11}|u|^2dxdt.
\end{split}
\end{equation}

\textbf{Estimate of $J_2$.}

In $\Omega_2$,
$$\Psi\leq \gamma R^{2/3}(\frac{3}{\sqrt{T}})^{4/3}+2R^2\leq 5\gamma(\frac{R}{T})^{2/3}+2R^2,\\$$
and by (\ref{pt2a}),
$$|\tilde{u}|+|\nabla \tilde{u}|\leq C(n,\Lambda,\lambda,M,N)e^{18N}\leq C,$$
thus
\begin{equation}\label{rhs1}
J_2\leq C T^{-\frac{n}{2}}\exp\{10\gamma(\frac{R}{T})^{2/3}+4R^2\}.
\end{equation}

\textbf{Estimate of $J_3$.}
In $\Omega_3$, $\psi(t)=0$,
$$\Psi\leq\gamma R^{2/3}(\gamma^{-3/4}R)^{4/3}=R^2,$$
and by (\ref{pt2a}),
$$|\tilde{u}|+|\nabla \tilde{u}|\leq C\exp\{2NT(\gamma^{-3/4}R)^2\}\leq C\exp\{2N(\gamma^{-3/4}R)^2\} ,$$
then we have
$$J_3\leq C(\gamma^{-3/4}R)^n\exp\{2R^2+4N(\gamma^{-3/4}R)^2\}.$$
Notice that if $\gamma^{-3/4}R>C(n,N)$, then
$$(\gamma^{-3/4}R)^n\leq \exp\{N(\gamma^{-3/4}R)^2\},$$
 hence
\begin{equation}\label{rhs2}
J_3\leq C\exp\{2R^2+5N(\gamma^{-3/4}R)^2\}\leq C\exp\{3R^2\}.
\end{equation}

\textbf{Estimate of $J_4$.}

In $\Omega_4$,
$$\Psi\leq\gamma R^{2/3}(\gamma^{-3/4}R)^{4/3}+2R^2=3R^2,$$
then
\begin{equation}\label{rhs3}
\begin{split}
J_4\leq &\exp\{6R^2\}\int_{\Omega_4}(|\tilde{u}|^2+|\nabla \tilde{u}|^2)dxdt\\
\leq &T^{-n/2-1}\exp\{6R^2\}\int_{T/8}^{7T/8}
\int_{\gamma^{-3/4}\sqrt{T}R-1<|x|<\gamma^{-3/4}\sqrt{T}R}(|u|^2+|\nabla u|^2)dxdt.\\
\end{split}
\end{equation}

Now we combine (\ref{ideab}), (\ref{lhsa}), (\ref{rhs1}), (\ref{rhs2}) and (\ref{rhs3})
together, then we have
\begin{eqnarray*}
\begin{split}
&\exp\{12\gamma (\frac{R}{T})^{2/3}+4R^2\}\frac{1}{T}\int_{T/3}^{2T/3}\int_{9\leq|x|\leq11}|u|^2dxdt\\
\leq &C\gamma^{-1}[\exp\{10\gamma(\frac{R}{T})^{2/3}+4R^2\}\\
&~~~~~~~+\exp\{6R^2\}\frac{1}{T}\int_{T/8}^{7T/8}
\int_{\gamma^{-3/4}\sqrt{T}R-1<|x|<\gamma^{-3/4}\sqrt{T}R}(|u|^2+|\nabla u|^2)dxdt].
\end{split}
\end{eqnarray*}
In the above inequality, we divide both sides by $\exp\{10\gamma (\frac{R}{T})^{2/3}+4R^2\}$,
and take\\ $\gamma=\gamma(n,\Lambda,\lambda,M,E,N)$ large enough,
then when $\gamma^{-3/4}R\geq\frac{C(n,N)}{\sqrt{T}}$, we have
\begin{eqnarray*}
\begin{split}
&\exp\{2\gamma(\frac{R}{T})^{2/3}\}\frac{1}{T}\int_{T/3}^{2T/3}\int_{9\leq|x|\leq11}|u|^2dxdt\\
\leq&1+\exp\{2R^2\}\frac{1}{T}\int_{T/8}^{7T/8}
\int_{\gamma^{-3/4}\sqrt{T}R-1<|x|<\gamma^{-3/4}\sqrt{T}R}(|u|^2+|\nabla u|^2)dxdt,
\end{split}
\end{eqnarray*}
If we replace $\gamma^{-3/4}\sqrt{T}R$ by $R$, we rewrite the above formula as
\begin{eqnarray*}
\begin{split}
&\exp\{2\gamma^{3/2}\frac{R^{2/3}}{T}\}\frac{1}{T}\int_{T/3}^{2T/3}\int_{9\leq|x|\leq11}|u|^2dxdt\\
\leq&1+\exp\{2\gamma^{3/2}\frac{R^2}{T}\}
\frac{1}{T}\int_{T/8}^{7T/8}\int_{R-1<|x|<R}(|u|^2+|\nabla u|^2)dxdt,
\end{split}
\end{eqnarray*}
provided $\gamma=\gamma(n,\Lambda,\lambda,M,E,N)$ large enough, and $R\geq C(n,N)$.
Thus we proved Lemma \ref{glbound}.\\

\emph{Proof of Proposition \ref{lbound}}. Since $u(x,0)\neq0$, then
by the unique continuation (see \cite{EF,F}), we must have
$u(x,0)\neq0$ in $B(10e_1,\frac{1}{2})$, and thus
$||u(\cdot,0)||_{L^2(B(10e_1,\frac{1}{2}))}>0.$

Now we apply Lemma \ref{doubleproperty} for $\rho=\frac{1}{2}$ and
the ball $B(10e_1,\frac{1}{2})$, then when
$$0<t\leq1/C\log(C\Theta_{1/2}),$$ we have
\begin{equation}\label{double}
C\int_{B(10e_1,1)}u^2(x,t)dx\geq \int_{B(10e_1,\frac{1}{2})}u^2(x,0)dx.
\end{equation}
Notice that $$\Theta_{1/2}\leq
C(N)/||u(\cdot,0)||^2_{L^2(B(10e_1,\frac{1}{2}))},$$ and
$$1/C\log(C\Theta_{1/2})\geq C(n,\Lambda,\lambda,M,N,||u(\cdot,0)||_{L^2(B(10e_1,\frac{1}{2}))})
\equiv T_2,$$ then when $0<t\leq T_2$, we have (\ref{double}).

For $0<T\leq T_2$, we apply Lemma \ref{glbound}, then when $R\geq R_3(n,N)$, we have
\begin{equation}\label{lowerbounda}
\begin{split}
&e^{C_\star\frac{R^{2/3}}{T}}\frac{1}{T}\int_{T/3}^{2T/3}\int_{9\leq|x|\leq11}u^2dxdt\\
\leq&1+e^{C_\star\frac{R^2}{T}}\frac{1}{T}\int_{T/8}^{7T/8}\int_{R-1\leq|x|\leq R}(u^2+|\nabla u|^2)dxdt,
\end{split}
\end{equation}
Notice that the left hand side of (\ref{lowerbounda})
$$e^{C_\star\frac{R^{2/3}}{T}}\frac{1}{T}\int_{T/3}^{2T/3}\int_{9\leq|x|\leq11}u^2dxdt
\geq e^{C_\star\frac{R^{2/3}}{T}}\frac{1}{T}\int_{T/3}^{2T/3}\int_{B(10e_1,1)}u^2dxdt,$$
and by (\ref{double}),
$$e^{C_\star\frac{R^{2/3}}{T}}\frac{1}{T}\int_{T/3}^{2T/3}\int_{9\leq|x|\leq11}u^2dxdt
\geq
Ce^{C_\star\frac{R^{2/3}}{T}}||u(\cdot,0)||^2_{L^2(B(10e_1,\frac{1}{2}))}.$$
If we choose $$R\geq
R_2(n,\Lambda,\lambda,M,E,N,||u(\cdot,0)||_{L^2(B(10e_1,\frac{1}{2}))}),$$
 then
\begin{equation}\label{lhsb}
e^{C_\star\frac{R^{2/3}}{T}}\frac{1}{T}\int_{T/3}^{2T/3}\int_{9\leq|x|\leq11}u^2dxdt
\geq e^{C_\star\frac{R^{2/3}}{2T}}\geq2.
\end{equation}
By (\ref{lowerbounda}) and (\ref{lhsb}), we have
$$\frac{1}{T}\int_{T/8}^{7T/8}\int_{R-1\leq|x|\leq R}(u^2+|\nabla u|^2)dxdt
\geq e^{-C_\star\frac{R^2}{T}}.$$ Thus we proved the lower bound
estimate.

\section{Proof of Carleman Inequalities}

In this section, we shall prove the two Carleman Inequalities. The main idea is to choose a
proper weighted functions $G$. We denote
$$\tilde{\Delta}v=\partial_i(a^{ij}\partial_jv).$$ Here and in the following argument,
we use the summation convention on the repeated indices. We shall
make use of the following lemma which is due to Escauriaza and
Fern\'{a}ndez in \cite{EF} (see also \cite{Tu}).
\begin{Lemma}\label{lem-generalCI}
Suppose $\sigma(t):\mathbb{R}_+\rightarrow \mathbb{R}_+$ is a smooth
function, $F$ is differentiable, $G$ is twice differentiable and $G>0$.
Then the following identity holds for any $v\in
C^\infty_0(\mathbb{R}^n\times[0,T])$ and any $\alpha\in\mathbb{R}$:
\begin{equation}\label{generalCI}
\begin{split}
&2\int_{\mathbb{R}^n\times[0,T]}\frac{\sigma^{1-\alpha}}{\sigma'}|Lv|^2Gdxdt
+\frac{1}{2}\int_{\mathbb{R}^n\times[0,T]}\frac{\sigma^{1-\alpha}}{\sigma'}v^2MGdxdt\\
&+\int_{\mathbb{R}^n\times[0,T]}\frac{\sigma^{1-\alpha}}{\sigma'}
\langle \mathbf{A}\nabla v,\nabla v\rangle
[(\log{\frac{\sigma}{\sigma'}})'+\frac{\partial_tG-\tilde{\Delta}G}{G}-F]Gdxdt\\
&+2\int_{\mathbb{R}^n\times[0,T]}\frac{\sigma^{1-\alpha}}{\sigma'}
\langle \mathbf{D}_G\nabla v,\nabla v\rangle Gdxdt
-\int_{\mathbb{R}^n\times[0,T]}\frac{\sigma^{1-\alpha}}{\sigma'}v
\langle \mathbf{A}\nabla v,\nabla F\rangle Gdxdt\\
=&2\int_{\mathbb{R}^n\times[0,T]}\frac{\sigma^{1-\alpha}}{\sigma'}LvPvGdxdt
+\int_{\mathbb{R}^n}\frac{\sigma^{1-\alpha}}{\sigma'}\langle \mathbf{A}\nabla v,\nabla v\rangle Gdx|^T_0\\
&~~~~~~~+\frac{1}{2}\int_{\mathbb{R}^n}\frac{\sigma^{1-\alpha}}{\sigma'}v^2(F-
\frac{\alpha\sigma'}{\sigma})Gdx|^T_0
\end{split}
\end{equation}
where
$$Lv=\partial_tv-\langle \mathbf{A}\nabla{logG},\nabla v\rangle+\frac{1}{2}
(F-\frac{\alpha\sigma'}{\sigma})v,$$
$$M=(log{\frac{\sigma}{\sigma'}})'F+\partial_tF+(F-\frac{\alpha\sigma'}{\sigma})(\frac{\partial_tG-
\tilde{\Delta}G}{G}-F)-\langle \mathbf{A}\nabla F,\nabla{logG}\rangle,$$
and
$$\mathbf{D}^{ij}_G=a^{ik}\partial_{kl}(logG)a^{lj}
+\frac{\partial_l(logG)}{2}(a^{ki}\partial_ka^{lj}+a^{kj}\partial_ka^{li}-a^{kl}\partial_ka^{ij})+
\frac{1}{2}\partial_ta^{ij}.
$$
\end{Lemma}

We first give a modification of this lemma which will be used in our
proof. Let $\alpha=0$ and $\sigma(t)=e^t$ in Lemma
\ref{lem-generalCI}, then we obtain the following identity for $v\in
C^\infty_0(\mathbb{R}^n\times[0,T])$
\begin{eqnarray*}
\begin{split}
&\frac{1}{2}\int_{\mathbb{R}^n\times[0,T]}v^2MGdxdt
+\int_{\mathbb{R}^n\times[0,T]}\langle[2\mathbf{D}_G+
(\frac{\partial_tG-\tilde{\Delta}G}{G}-F)\mathbf{A}]\nabla v,\nabla v\rangle Gdxdt\\
&-\int_{\mathbb{R}^n\times[0,T]}v\langle \mathbf{A}\nabla v,\nabla F\rangle Gdxdt=
2\int_{\mathbb{R}^n\times(0,T)}Lv(Pv-Lv)Gdxdt\\
&~~~~~~~~~~~~~~~~~~~~~~~~~~+
\int_{\mathbb{R}^n}\langle \mathbf{A}\nabla v,\nabla v\rangle Gdx|^T_0+
\frac{1}{2}\int_{\mathbb{R}^n}v^2FGdx|^T_0.
\end{split}
\end{eqnarray*}
If $\nabla F$ is differentiable, we can integrate by parts to obtain that
\begin{eqnarray*}
\begin{split}
&-\int_{\mathbb{R}^n\times[0,T]}v\langle \mathbf{A}\nabla v,\nabla F\rangle Gdxdt\\
=&\frac{1}{2}\int_{\mathbb{R}^n\times[0,T]}v^2\tilde{\Delta} FGdxdt+
\frac{1}{2}\int_{\mathbb{R}^n\times[0,T]}v^2\langle \mathbf{A}\nabla F,\nabla logG\rangle Gdxdt.
\end{split}
\end{eqnarray*}
The function $\nabla F$ may not be differentiable, so we approximate
$F$ by some twice differentiable function $F_0$ and use the above
identity with $F_0$ in place of $F$, following Nguyen's idea in
\cite{Tu}.  Thus a direct corollary follows.

\begin{Corollary}\label{cor-generalCI1}
Suppose $F$ is differentiable, $F_0$ and $G$ is twice differentiable and $G>0$. Then the
following identity holds for any $v\in C^\infty_0(\mathbb{R}^n\times[0,T])$:
\begin{equation}\label{generalCI1}
\begin{split}
&\frac{1}{2}\int_{\mathbb{R}^n\times[0,T]}v^2M_0Gdxdt+
\int_{\mathbb{R}^n\times[0,T]}\langle [2\mathbf{D}_G+(\frac{\partial_tG-
\tilde{\Delta}G}{G}-F)\mathbf{A}]\nabla v,\nabla v\rangle Gdxdt\\
&-\int_{\mathbb{R}^n\times[0,T]}v\langle \mathbf{A}\nabla v,\nabla (F-F_0)\rangle Gdxdt=
2\int_{\mathbb{R}^n\times[0,T]}Lv(Pv-Lv)Gdxdt\\
&~~~~~~~~~~~~~~~+\int_{\mathbb{R}^n}\langle \mathbf{A}\nabla v,\nabla v\rangle Gdx|^T_0+
\frac{1}{2}\int_{\mathbb{R}^n}v^2FGdx|^T_0,
\end{split}
\end{equation}
where
$$Lu=\partial_tv-\langle \mathbf{A}\nabla v,\nabla{logG}\rangle+\frac{Fv}{2},$$
$$M_0=\partial_tF+F(\frac{\partial_tG-\tilde{\Delta}G}{G}-F)+\tilde{\Delta}F_0-
\langle \mathbf{A}\nabla (F-F_0),\nabla{logG}\rangle,$$
and
$$\mathbf{D}^{ij}_G=a^{ik}\partial_{kl}(logG)a^{lj}
+\frac{\partial_l(logG)}{2}(a^{ki}\partial_ka^{lj}+a^{kj}\partial_ka^{li}-
a^{kl}\partial_ka^{ij})+\frac{1}{2}\partial_ta^{ij}.
$$
\end{Corollary}

Before we prove our Carleman inequalities, we need to prove a result which can be viewed as another
version of Corollary \ref{cor-generalCI1}.\\

In (\ref{generalCI1}), we let $G=e^{2\Phi}$, $w=e^{\Phi}v$, and we denote
$$\mathbf{B}=2\mathbf{D}_G+(\frac{\partial_tG-\tilde{\Delta}G}{G}-F)\mathbf{A}.$$
Then the third term of the left hand side of (\ref{generalCI1}) is
\begin{eqnarray*}
\begin{split}
&-\int_Qv\langle \mathbf{A}\nabla(F-F_0),\nabla v\rangle e^{2\Phi}dxdt\\
=&-\int_Qw\langle \mathbf{A}\nabla(F-F_0),\nabla w-\nabla \Phi w\rangle dxdt\\
=&-\int_Qw\langle \mathbf{A}\nabla(F-F_0),\nabla w\rangle dxdt+
\int_Q\langle \mathbf{A}\nabla(F-F_0),\nabla \Phi \rangle w^2dxdt.\\
\end{split}
\end{eqnarray*}
We use the above identity and rewrite (\ref{generalCI1}) as
\begin{equation}\label{generalCI2}
\begin{split}
&\frac{1}{2}\int_QM_1w^2dxdt+\int_Q\langle \mathbf{B}\nabla v,\nabla v\rangle e^{2\Phi}dxdt-
\int_Qw\langle \mathbf{A}\nabla(F-F_0),\nabla w\rangle dxdt\\
=&2\int_QLv(Pv-Lv)e^{2\Phi}dxdt+\int_{\mathbb{R}^n}
\langle \mathbf{A}\nabla v,\nabla v\rangle e^{2\Phi}dx|^T_0+\frac{1}{2}
\int_{\mathbb{R}^n}v^2Fe^{2\Phi}dx|^T_0
\end{split}
\end{equation}
where
$$M_1=\partial_tF+F(\frac{\partial_tG-\tilde{\Delta}G}{G}-F)+\tilde{\Delta}F_0,$$
\begin{equation}\label{B}
\mathbf{B}=4\mathbf{A}D^2\Phi \mathbf{A}+2\partial_l\Phi(a^{ki}\partial_ka^{lj}+
a^{kj}\partial_ka^{li}-a^{kl}\partial_ka^{ij})
+\partial_ta^{ij}+(\frac{\partial_tG-\tilde{\Delta}G}{G}-F)\mathbf{A}.
\end{equation}
By direct calculations we have
\begin{equation}\label{G}
\frac{\partial_tG-\tilde{\Delta}G}{G}=2\partial_t\Phi-2a^{ij}\partial_{ij}\Phi-
2\partial_ia^{ij}\partial_j\Phi-4\langle \mathbf{A}\nabla\Phi,\nabla\Phi\rangle.
\end{equation}
Let
\begin{equation}\label{choiceofF}
F=2\partial_t\Phi-2a^{ij}\partial_{ij}\Phi-4\langle \mathbf{A}\nabla\Phi,\nabla\Phi\rangle-H,
\end{equation}
where $H$ is a smooth function to be determined. We choose
$$F_0=2\partial_t\Phi-2a^{ij}_\epsilon\partial_{ij}\Phi-4a^{ij}_\epsilon\partial_i\Phi\partial_j\Phi-H,$$
where
$$a^{ij}_\epsilon(x,t)=\int_{\mathbb{R}^n}a^{ij}(x-y,t)\phi_\epsilon(y)dy,$$
$\phi$ is a mollifier, and $\epsilon=\frac{1}{2}$.

By (\ref{B})-(\ref{choiceofF}), we have
\begin{equation}\label{B1}
\mathbf{B}=4\mathbf{A}D^2\Phi \mathbf{A}+2\partial_l\Phi(a^{ki}\partial_ka^{lj}+
a^{kj}\partial_ka^{li}-a^{kl}\partial_ka^{ij}-a^{ij}\partial_ka^{kl})
+\partial_ta^{ij}+H\mathbf{A}.
\end{equation}
Now we begin to prove our Carleman inequalities.

\subsection{Proof of Proposition \ref{Prop-C1}.}

Note that Carleman inequality (\ref{C1}) is very similar to the
second Carleman inequality in \cite{WZ}, and their proofs are also
similar.

In this part, we let
$$\Phi=\gamma f(t)|x|^{3/2}-\frac{b|x|^2+\beta}{2(t+1)},$$
where $b=\frac{1}{16\Lambda}$ and $\beta=\beta(n,\Lambda,\lambda,M,E)$ large enough.\\

\textbf{Step 1.} Estimate matrix $\mathbf{B}$.

We need to estimate the lower bounds of the matrices in the right side of (\ref{B1}).

First we estimate $D^2\Phi$. Denote that $$h=\gamma f|x|^{-1/2}.$$
By direct calculations we have
\begin{eqnarray*}
\begin{split}
D^2\Phi=\frac{3}{2}h (I_n-\frac{x\cdot x^T}{2|x|^2})-\frac{b}{t+1}I_n\geq(\frac{3}{4}h-\frac{b}{t+1})I_n,
\end{split}
\end{eqnarray*}
and hence
$$4\textbf{A}D^2\Phi \textbf{A}\geq(3\lambda^2h-\frac{C}{t+1})I_n.$$

Second, we estimate matrix $\partial_l\Phi a^{ki}\partial_ka^{lj}$
and $\partial_ta^{ij}$. For any $\xi\in\mathbb{R}^n$,
$$
|\partial_l\Phi a^{ki}\partial_ka^{lj}\xi_i\xi_j|
\leq n^2\Lambda\frac{E}{|x|}|\nabla\Phi|\sum_{i,j}|\xi_i||\xi_j|
\leq \frac{n^3\Lambda E}{|x|}|\nabla\Phi||\xi|^2,
$$
Since
\begin{equation}\label{pCI1-graphi}
\nabla\Phi=(\frac{3}{2}h -\frac{b}{t+1})x,
\end{equation}
then
$$
|\partial_l\Phi a^{ki}\partial_ka^{lj}\xi_i\xi_j|\leq n^{3}\Lambda E(\frac{3}{2}h+\frac{b}{t+1}),
$$
and thus
$$\partial_l\Phi a^{ki}\partial_ka^{lj}\geq -n^{3}\Lambda E(\frac{3}{2}h+\frac{b}{t+1})I_n.$$
Similarly,
$$\partial_ta^{ij}\geq-nMI_n.$$
Consequently,
\begin{eqnarray*}
\begin{split}
\textbf{B}\geq&(3\lambda^2-12n^3\Lambda E)h I_n-\frac{C}{t+1} I_n-n MI_n+\lambda HI_n\\
\geq&2\lambda^2h I_n+(\lambda H-\frac{C}{t+1})I_n,
\end{split}
\end{eqnarray*}
if we take $E<E_0(n,\Lambda,\lambda)$. Now in this part, we choose
$$H=\frac{d}{t+1} ,$$
where $d=d(n,\Lambda,\lambda,M,E)$ large enough, then we have
\begin{equation}\label{pCI1-estB}
\textbf{B}\geq2\lambda^2(h+\frac{1}{t+1})I_n+I_n.
\end{equation}

\textbf{Step 2.} Prove the Carleman inequality.

By (\ref{pCI1-estB}), we can estimate the second term of the left hand side of (\ref{generalCI2}),
\begin{equation}\label{pCI1a}
\begin{split}
&\int_{Q}\langle \mathbf{B}\nabla v,\nabla v\rangle e^{2\Phi}dxdt\\
\geq&\int_{Q}e^{2\Phi}|\nabla v|^2dxdt
+2\lambda^2\int_{Q}(h+\frac{1}{t+1})e^{2\Phi}|\nabla v|^2dxdt\\
=&\int_{Q}e^{2\Phi}|\nabla v|^2dxdt+2\lambda^2\int_{Q}(h+\frac{1}{t+1})|\nabla w|^2dxdt\\
&+2\lambda^2\int_{Q}[(h+\frac{1}{t+1})|\nabla\Phi|^2+\nabla h\cdot\nabla\Phi+(h+
\frac{1}{t+1})\Delta\Phi]w^2dxdt.
\end{split}
\end{equation}
By (\ref{generalCI2}), (\ref{pCI1a}) and the Cauchy inequality, we have
\begin{equation}\label{pCI1b}
\begin{split}
&\int_{Q}e^{2\Phi}|\nabla v|^2dxdt+2\lambda^2\int_{Q}(h+\frac{1}{t+1})|\nabla w|^2dxdt+
\int_{Q}M_2w^2dxdt\\
&-\int_{Q}w\langle A\nabla(F-F_0),\nabla w\rangle dxdt\leq \int_{Q}e^{2\Phi}|Pv|^2dxdt\\
&~~~~~~~~~~+\int_{\mathbb{R}^n}\langle \mathbf{A}\nabla v,\nabla v\rangle e^{2\Phi}dx|^1_0+
\frac{1}{2}\int_{\mathbb{R}^n}Fe^{2\Phi}v^2dx|^1_0,
\end{split}
\end{equation}
where
\begin{equation}
\begin{split}
M_2=&2\lambda^2[(h+\frac{1}{t+1})|\nabla\Phi|^2+\nabla h\cdot\nabla\Phi+(h+\frac{1}{t+1})\Delta\Phi]\\
&+\frac{1}{2}\partial_tF+\frac{1}{2}F(\frac{\partial_tG-\tilde{\Delta}G}{G}-F)+
\frac{1}{2}\tilde{\Delta}F_0.
\end{split}
\end{equation}

We use inequality (\ref{pCI1b}) to prove Proposition \ref{Prop-C1}.
We need some estimates which we list in the following lemma.

\begin{Lemma}\label{pCI1-estimates}
Set $b=\frac{1}{16\Lambda}$, $\beta=20\frac{\Lambda}{\lambda}d$ and $d=d(n,\Lambda,\lambda,M,E)$
large enough. There exists $E_0=E_0(n,\Lambda,\lambda)$,
such that when $E<E_0$, for any $\gamma>0$, we have
\begin{equation}\label{pCI1-estM2}
M_2\geq \lambda^2h^3|x|^2+\frac{db}{8}\frac{|x|^2}{(t+1)^3};
\end{equation}
\begin{equation}\label{pCI1-esterror}
|\nabla(F-F_0)|\leq C(n)E[h^2+\frac{1}{(t+1)^2}]|x|;
\end{equation}
\begin{equation}\label{pCI1-estF0}
F(x,0)\geq-2\beta|x|^{2}(1+\gamma)^2;
\end{equation}
\begin{equation}\label{pCI1-estF1}
F(x,1)\leq\frac{\beta}{2}|x|^2.
\end{equation}
\end{Lemma}

We shall prove this lemma later.

By applying Lemma 3.3, in particular by (\ref{pCI1-esterror}), we
have
\begin{eqnarray*}
\begin{split}
|\int_{Q}w\langle A\nabla(F-F_0),\nabla w\rangle dx dt|
\leq&\Lambda\int_{Q}|\nabla(F-F_0)||w||\nabla w|dxdt\\
\leq&C(n)\Lambda E\int_{Q}[h^2+\frac{1}{(t+1)^2}]|x||w||\nabla w|dx dt.
\end{split}
\end{eqnarray*}
Using the Cauchy inequality, we have
\begin{eqnarray*}
\begin{split}
|\int_{Q}w\langle A\nabla(F-F_0),\nabla w\rangle dx dt|\leq& C(n)\Lambda E\int_{Q}(h^3|x|^2+
\frac{|x|^2}{(t+1)^3})w^2dxdt\\
&+C(n)\Lambda E\int_{Q}(h+\frac{1}{t+1})|\nabla w|^2dxdt.
\end{split}
\end{eqnarray*}
When $E<E_0(n,\Lambda,\lambda)$, we have
\begin{equation}\label{pCI1c}
\begin{split}
|\int_{Q}w\langle A\nabla(F-F_0),\nabla w\rangle dx dt|\leq& \lambda^2\int_{Q}(h^3|x|^2+
\frac{|x|^2}{(t+1)^3})w^2dxdt\\
&+\lambda^2\int_{Q}(h+\frac{1}{t+1})|\nabla w|^2dxdt.
\end{split}
\end{equation}
Because of (\ref{pCI1b}), (\ref{pCI1c}) and (\ref{pCI1-estM2}), we have
\begin{equation}
\begin{split}\label{pCI1d}
&\int_{Q}e^{2\Phi}|\nabla v|^2dxdt+(\frac{db}{8}-C)\int_{Q}\frac{|x|^2}{(t+1)^3}w^2dxdt
\leq \int_{Q}e^{2\Phi}|Pv|^2dxdt\\
&~~~~~~~~~~+\int_{\mathbb{R}^n}\langle \mathbf{A}\nabla v,\nabla v\rangle e^{2\Phi}dx|^1_0+
\frac{1}{2}\int_{\mathbb{R}^n}Fe^{2\Phi}v^2dx|^1_0.
\end{split}
\end{equation}
Now we estimate the second term of the right hand side of (\ref{pCI1d}).
\begin{equation*}
\begin{split}
\int_{\mathbb{R}^n}\langle \mathbf{A}\nabla v,\nabla v\rangle e^{2\Phi}dx|^1_0
\leq&\int_{\mathbb{R}^n}\langle \mathbf{A}\nabla v,\nabla v\rangle e^{2\Phi}dx|_{t=1}\\
\leq&\Lambda\int_{\mathbb{R}^n}|\nabla v(x,1)|^2e^{-\frac{b|x|^2+\beta}{2}}dxdt
\end{split}
\end{equation*}
Notice that
$$\Lambda e^{-\frac{\beta}{2}}\leq \beta e^{-\frac{\beta}{2}}\leq c,$$
then
\begin{equation}\label{pCI1e}
\int_{\mathbb{R}^n}\langle \mathbf{A}\nabla v,\nabla v\rangle e^{2\Phi}dx|^1_0
\leq c\int_{\mathbb{R}^n}|\nabla v(x,1)|^2e^{-\frac{b|x|^2}{2}}dxdt
\end{equation}
Finally we estimate the third term of the right hand side of
(\ref{pCI1d}).
\begin{eqnarray*}
\begin{split}
\frac{1}{2}\int_{\mathbb{R}^n}Fe^{2\Phi}v^2dx|^1_0=&\frac{1}{2}
\int_{\mathbb{R}^n}F(x,1)e^{-\frac{b|x|^2+\beta}{2}}v^2(x,1)dx\\
&-\frac{1}{2}\int_{\mathbb{R}^n}F(x,0)e^{2\gamma (1-2^{-\beta})|x|^{3/2}-b|x|^2-\beta}v^2(x,0)dx
\end{split}
\end{eqnarray*}
By (\ref{pCI1-estF0}) and (\ref{pCI1-estF1}), we have
\begin{equation}\label{pCI1f}
\begin{split}
\frac{1}{2}\int_{\mathbb{R}^n}Fe^{2\Phi}v^2dx|^1_0
\leq&\frac{\beta}{4}\int_{\mathbb{R}^n}|x|^2e^{-\frac{b|x|^2+\beta}{2}}v^2(x,1)dx\\
&+\beta(1+\gamma)^2\int_{\mathbb{R}^n}|x|^2e^{2\gamma (1-2^{-\beta})|x|^{3/2}-b|x|^2-\beta}v^2(x,0)dx\\
\leq&c\int_{\mathbb{R}^n}|x|^2e^{-\frac{b|x|^2}{2}}v^2(x,1)dx\\
&+c(1+\gamma)^2\int_{\mathbb{R}^n}|x|^2e^{2\gamma|x|^{3/2}-b|x|^2}v^2(x,0)dx.
\end{split}
\end{equation}
We combine (\ref{pCI1d}), (\ref{pCI1e}) and (\ref{pCI1f}), and take
$d$ large enough, then we proved Carleman inequality (\ref{C1}).

It remains to prove Lemma \ref{pCI1-estimates}.\\

\textbf{Step 3.} Prove Lemma \ref{pCI1-estimates}.

\textbf{Estimate of $M_2$.}

We estimate the terms of $M_2$ respectively.\\

\emph{Estimate of the first three terms.}

By (\ref{pCI1-graphi}), we have
\begin{eqnarray*}
\begin{split}
(h+\frac{1}{t+1})|\nabla\Phi|^2\geq &h|\nabla\Phi|^2=h(\frac{3}{2}h-\frac{b}{t+1})^2|x|^2\\
\geq &h|x|^2[\frac{9}{8}h^2-\frac{b^2}{(t+1)^2}]\\
=&[\frac{9}{8}h^3-\frac{C}{(t+1)^2}h]|x|^2;
\end{split}
\end{eqnarray*}

$$\nabla h\cdot\nabla\Phi=-\frac{1}{2}h(\frac{3}{2}h-\frac{b}{t+1})\geq-\frac{3}{4}h^2;$$

$$(h+\frac{1}{t+1})\Delta\Phi\geq-\frac{nb}{t+1}(h+\frac{1}{t+1})
\geq-\frac{C}{t+1}h-\frac{C}{(t+1)^2}.$$ Combining them together, we
obtain
\begin{equation}\label{aaa1}
2\lambda^2(h|\nabla\Phi|^2+\nabla h\cdot\nabla\Phi+h\Delta\Phi)
\geq [\frac{9}{4}\lambda^2 h^3-Ch^2-\frac{C}{(t+1)^2}h]|x|^2-\frac{C}{(t+1)^3}.
\end{equation}

\emph{Estimate of $\frac{1}{2}\partial_tF$.}

Recall (\ref{choiceofF}), then
$$\frac{1}{2}\partial_tF=\partial_t^2\Phi-\partial_ta^{ij}\partial_{ij}\Phi-a^{ij}\partial_{ijt}\Phi
-2\partial_t\langle A\nabla\Phi,\nabla\Phi\rangle-\frac{1}{2}\partial_t H.$$
We estimate them one by one. Keep in mind that $f'<0$.
$$\partial_t^2\Phi=\gamma f''|x|^{3/2}-\frac{b|x|^2+\beta}{(t+1)^3}=
\frac{f''}{f}h|x|^2-\frac{b|x|^2+\beta}{(t+1)^3};$$

$$
-\partial_ta^{ij}\partial_{ij}\Phi=-\frac{3}{2}h(\partial_ta^{ii}-
\frac{\partial_ta^{ij}x_ix_j}{2|x|^2})+\frac{b\partial_ta^{ii}}{t+1}
\geq-C h-\frac{C}{t+1};
$$

$$
-a^{ij}\partial_{ijt}\Phi=-\frac{3f'}{2f}h(a^{ii}-\frac{a^{ij}x_ix_j}{2|x|^2})-\frac{ba^{ii}}{(t+1)^2}
\geq C \frac{f'}{f}h-\frac{C}{(t+1)^2};
$$
\begin{eqnarray*}
\begin{split}
&-2\partial_t\langle A\nabla\Phi,\nabla\Phi\rangle\\
=&[-9\frac{f'}{f}h^2+6b(\frac{f'}{(t+1)f}-\frac{1}{(t+1)^2})h+\frac{4b^2}{(t+1)^3}]a^{ij}x_ix_j\\
&-2(\frac{3}{2}h-\frac{b}{t+1})^2\partial_ta^{ij}x_ix_j\\
\geq&[-9\lambda\frac{f'}{f}h^2+C(\frac{f'}{(t+1)f}-\frac{1}{(t+1)^2})h]|x|^2-
C[h^2+\frac{1}{(t+1)^2}]|x|^2\\
\geq&[(-9\lambda\frac{f'}{f}-C)h^2+C(\frac{f'}{(t+1)f}-\frac{1}{(t+1)^2})h]|x|^2-\frac{C|x|^2}{(t+1)^2};
\end{split}
\end{eqnarray*}
$$-\frac{1}{2}\partial_t H=\frac{d}{2(t+1)^2}.$$
Combining them together, we have
\begin{equation}\label{aaa2}
\frac{1}{2}\partial_tF\geq[(-9\lambda\frac{ f'}{f}-C)h^2+
(\frac{f''}{f}+\frac{Cf'}{(t+1)f}-\frac{C}{(t+1)^2})h]|x|^2
-\frac{C|x|^2+\beta}{(t+1)^3};
\end{equation}

\emph{Estimate of $\frac{1}{2}F(\frac{\partial_tG-\tilde{\Delta}G}{G}-F)$.}

First we have
\begin{eqnarray*}
\begin{split}
\frac{1}{2}F(\frac{\partial_tG-\tilde{\Delta}G}{G}-F)
=&(\partial_t\Phi-2\langle A\nabla\Phi,\nabla\Phi\rangle-a^{ij}\partial_{ij}\Phi-
\frac{1}{2}H)(H-2\partial_ia^{ij}\partial_j\Phi)\\
\equiv&J_1-J_2-J_3,
\end{split}
\end{eqnarray*}
where
\begin{eqnarray*}
\begin{split}
J_1=&\partial_t\Phi(H-2\partial_ia^{ij}\partial_j\Phi)\\
J_2=&2\langle A\nabla\Phi,\nabla\Phi\rangle(H-2\partial_ia^{ij}\partial_j\Phi)\\
J_3=&(a^{ij}\partial_{ij}\Phi+\frac{1}{2}H)(H-2\partial_ia^{ij}\partial_j\Phi).
\end{split}
\end{eqnarray*}
Before we estimate $J_1$, $J_2$ and $J_3$, we estimate $2\partial_ia^{ij}\partial_j\Phi$ first.
$$|2\partial_ia^{ij}\partial_j\Phi|\leq\frac{2n^2E}{|x|}|\nabla\Phi|
\leq2n^2E(\frac{3}{2}h+\frac{b}{t+1})\leq 3n^2Eh+\frac{C}{t+1},$$
then
$$-3n^2Eh+\frac{d-C}{t+1}\leq H-2\partial_ia^{ij}\partial_j\Phi\leq 3n^2Eh+\frac{d+C}{t+1}.$$
Now we estimate $J_1$, $J_2$ and $J_3$ respectively.
\begin{eqnarray*}
\begin{split}
J_1=&\frac{f'}{f}h|x|^2(H-2\partial_ia^{ij}\partial_j\Phi)+\frac{b|x|^2+
\beta}{2(t+1)^2}(H-2\partial_ia^{ij}\partial_j\Phi)\\
\geq&\frac{f'}{f}h|x|^2(3n^2Eh+\frac{d+C}{t+1})+\frac{b|x|^2+\beta}{2(t+1)^2}(-3n^2Eh+\frac{d-C}{t+1})\\
\geq&[3n^2E\frac{f'}{f}h^2+(\frac{(d+C) f'}{(t+1)f}-\frac{C\beta+C}{(t+1)^2})h]|x|^2+
(\frac{d}{2}-C)\frac{(b|x|^2+\beta)}{(t+1)^3},
\end{split}
\end{eqnarray*}
\begin{eqnarray*}
\begin{split}
J_2\leq&2\Lambda(\frac{3}{2}h-\frac{b}{t+1})^2|x|^2(3n^2Eh+\frac{d+C}{t+1})\\
\leq&4\Lambda[\frac{9}{4}h^2+\frac{b^2}{(t+1)^2}]|x|^2(3n^2Eh+\frac{d+C}{t+1})\\
\leq&[27n^2\Lambda Eh^3+\frac{C}{(t+1)^2} h+\frac{9d\Lambda +C}{t+1} h^2]|x|^2+
\frac{4d\Lambda b^2+C}{(t+1)^3}|x|^2,
\end{split}
\end{eqnarray*}
and
\begin{eqnarray*}
\begin{split}
J_3\leq&|\frac{3}{2}h(a^{ii}-\frac{a^{ij}x_ix_j}{2|x|^2})+\frac{d/2-
ba^{ii}}{t+1}|(3n^2Eh+\frac{d+C}{t+1})\\
\leq &(C h+\frac{d}{t+1})^2\leq Ch^2+\frac{2d^2}{(t+1)^2}\leq Ch^2|x|^2+\frac{4d^2}{(t+1)^3}.
\end{split}
\end{eqnarray*}
Combining $J_1$, $J_2$ and $J_3$ together, we obtain
\begin{equation}\label{aaa3}
\begin{split}
&\frac{1}{2}F(\frac{\partial_tG-\tilde{\Delta}G}{G}-F)\\
\geq& [-27n^2\Lambda Eh^3+(3n^2E\frac{f'}{f}-\frac{9\Lambda d+C}{t+1})h^2+
(\frac{(d+C)f'}{(t+1)f}-\frac{C\beta+C}{(t+1)^2})h]|x|^2\\
&+(\frac{db}{2}-4d\Lambda b^2-C)\frac{|x|^2}{(t+1)^3}+\frac{(d/2-C)\beta-4d^2}{(t+1)^3}.
\end{split}
\end{equation}

\emph{Estimate of $\frac{1}{2}\tilde{\Delta}F_0$.}

In order to estimate $\tilde{\Delta} F_0$ and $|\nabla(F-F_0)|$, we
need some estimates about $\{a^{ij}_\epsilon\}$ which we will prove in the appendix.

In fact, $\{a^{ij}_\epsilon\}$ satisfy the following properties:
\begin{equation}\label{indexA}
\begin{split}
&i)~\lambda|\xi|^2\leq a^{ij}_\epsilon(x,t)\xi_i\xi_j\leq \Lambda|\xi|^2, ~\forall \xi\in\mathbb{R}^n;\\
&ii)~|\nabla a^{ij}_\epsilon(x,t)|\leq M;~~|\nabla a^{ij}_\epsilon(x,t)|
\leq\frac{2E}{|x|}~when~|x|\geq1;\\
&iii)~|a^{ij}_\epsilon(x,t)-a^{ij}(x,t)|\leq 2\Lambda ;
~~|a^{ij}_\epsilon(x,t)-a^{ij}(x,t)|\leq\frac{E}{|x|}~when~|x|\geq1;\\
&iv)~|\partial_{kl}a^{ij}_\epsilon(x,t)|\leq c(n)M;~~|\partial_{kl}a^{ij}_\epsilon(x,t)|\leq
\frac{c(n)E}{|x|}~when~|x|\geq1.\\
\end{split}
\end{equation}
Recall that
\begin{eqnarray*}
\begin{split}
F_0=&2\partial_t\Phi-2a^{ij}_\epsilon\partial_{ij}\Phi-4a^{ij}_\epsilon\partial_i\Phi\partial_j\Phi-H\\\
=&2\gamma f'|x|^{3/2}+\frac{b|x|^2+\beta}{(t+1)^2}-\frac{d}{t+1}-
2a^{ij}_\epsilon(\partial_{ij}\Phi+2\partial_i\Phi\partial_j\Phi),
\end{split}
\end{eqnarray*}
then
\begin{equation}\label{aaa4a}
\begin{split}
\frac{1}{2}\tilde{\Delta}F_0=&\gamma f'\tilde{\Delta}(|x|^{3/2})+\frac{b}{2(t+1)^2}\tilde{\Delta}(|x|^2)
-\tilde{\Delta}[a^{ij}_\epsilon(\partial_{ij}\Phi+2\partial_i\Phi\partial_j\Phi)]\\
\geq&C\gamma f'|x|^{-1/2}-\frac{C}{(t+1)^2}-\tilde{\Delta}[a^{ij}_\epsilon(\partial_{ij}\Phi+
2\partial_i\Phi\partial_j\Phi)]\\
=&\frac{Cf'}{f}h-\frac{C}{(t+1)^2}-\tilde{\Delta}[a^{ij}_\epsilon(\partial_{ij}\Phi+
2\partial_i\Phi\partial_j\Phi)],
\end{split}
\end{equation}
and thus it remains to estimate
$\tilde{\Delta}[a^{ij}_\epsilon(\partial_{ij}\Phi+
2\partial_i\Phi\partial_j\Phi)]$.

By (\ref{indexA}) we have that $|a^{ij}_\epsilon|$, $|\nabla
a^{ij}_\epsilon|$ and $|\nabla^2a^{ij}_\epsilon|$ are all bounded,
and it is easy to verify that
\begin{equation}\label{aaa4b}
|\nabla^k\Phi|\leq C(n) (h+\frac{1}{t+1})|x|^{2-k},~~k=1,2,3,4.
\end{equation}
Direct calculations give us
\begin{eqnarray*}
\begin{split}
&|\tilde{\Delta}[a^{ij}_\epsilon(\partial_{ij}\Phi+2\partial_i\Phi\partial_j\Phi)]|\\
=&|\partial_ka^{kl}\partial_l[a^{ij}_\epsilon(\partial_{ij}\Phi+2\partial_i\Phi\partial_j\Phi)]
+a^{kl}\partial_{kl}[a^{ij}_\epsilon(\partial_{ij}\Phi+2\partial_i\Phi\partial_j\Phi)]|\\
\leq&C(|\nabla^2\Phi|+|\nabla^3\Phi|+|\nabla^4\Phi|+|\nabla\Phi|^2+
|\nabla\Phi||\nabla^2\Phi|+|\nabla\Phi||\nabla^3\Phi|+|\nabla^2\Phi|^2).\\
\end{split}
\end{eqnarray*}
By the Cauchy inequality, we have
\begin{equation}\label{estmollifier}
\begin{split}
&|\tilde{\Delta}[a^{ij}_\epsilon(\partial_{ij}\Phi+2\partial_i\Phi\partial_j\Phi)]|\\
\leq& C(|\nabla^2\Phi|+|\nabla^3\Phi|+|\nabla^4\Phi|+|\nabla\Phi|^2+|\nabla^2\Phi|^2+|\nabla^3\Phi|^2)
\end{split}
\end{equation}
then by (\ref{aaa4b}),
\begin{equation}\label{aaa4c}
|\tilde{\Delta}[a^{ij}_\epsilon(\partial_{ij}\Phi+2\partial_i\Phi\partial_j\Phi)]|
\leq C(h+\frac{1}{t+1})^2|x|^2
\leq C[h^2+\frac{1}{(t+1)^2}]|x|^2.
\end{equation}
We combine (\ref{aaa4a}) and (\ref{aaa4c}) and obtain that
\begin{equation}\label{aaa4}
\frac{1}{2}\tilde{\Delta}F_0\geq(\frac{Cf'}{f}h-Ch^2)|x|^2-\frac{C|x|^2}{(t+1)^3}.
\end{equation}

At last, combining (\ref{aaa1}), (\ref{aaa2}), (\ref{aaa3}) and (\ref{aaa4}),
we have
\begin{eqnarray*}
\begin{split}
M_2\geq &[(\frac{9}{4}\lambda^2-27n^2\Lambda E)h^3+(-(9\lambda-3n^2E)\frac{f'}{f}-
\frac{9\Lambda d+C}{t+1})h^2]|x|^2\\
&+(\frac{f''}{f}+\frac{(d+C)f'}{(t+1)f}-\frac{C\beta+C}{(t+1)^2})h|x|^2\\
&+(\frac{db}{2}-4d\Lambda b^2-C)\frac{|x|^2}{(t+1)^3}+\frac{(d/2-C)\beta-4d^2-C}{(t+1)^3}.\\
\end{split}
\end{eqnarray*}
Now we choose
$$b=\frac{1}{16\Lambda}, ~~~\beta=20\frac{\Lambda}{\lambda}d,$$
and we take $d$ large enough, then
$$\frac{db}{2}-4d\Lambda b^2-C=\frac{db}{4}-C\geq\frac{db}{8},$$
$$(\frac{d}{2}-C)\beta-4d^2-C\geq\frac{d}{4}\beta-5d^2\geq0, $$
and thus when $E<E_0(n,\Lambda,\lambda)$, we have
\begin{eqnarray*}
\begin{split}
M_2\geq &[\lambda^2h^3+(-\frac{9\lambda f'}{2f}-\frac{18\Lambda d}{t+1})h^2]|x|^2\\
&+(\frac{f''}{f}+\frac{2df'}{(t+1)f}-\frac{C\beta}{(t+1)^2})h|x|^2+\frac{db}{8}\frac{|x|^2}{(t+1)^3}.
\end{split}
\end{eqnarray*}
We take into account that
$$f(t)=(t+1)^{-\beta}-2^{-\beta},$$
then we have
\begin{eqnarray*}
\begin{split}
-\frac{9\lambda f'}{2f}-\frac{18\Lambda d}{t+1}=&\frac{9\lambda\beta}{2(t+1)[1-
(\frac{t+1}{2})^\beta]}-\frac{18\Lambda d}{t+1}\\
\geq&\frac{9\lambda\beta}{2(t+1)}-\frac{18\Lambda d}{t+1}=\frac{9(\lambda\beta-4\Lambda d)}{2(t+1)}\geq0,
\end{split}
\end{eqnarray*}
and
\begin{eqnarray*}
\begin{split}
\frac{f''}{f}+\frac{2df'}{(t+1)f}-\frac{C\beta}{(t+1)^2}=&\frac{\beta(\beta+1-2d)}{(t+1)^2[1-
(\frac{t+1}{2})^\beta]}-\frac{C\beta}{(t+1)^2}\\
\geq&\frac{\beta(\beta+1-2d)}{(t+1)^2}-\frac{d\beta}{(t+1)^2}\\
=&\frac{\beta(\beta+1-3d)}{(t+1)^2}\geq0,
\end{split}
\end{eqnarray*}
thus
\begin{equation}
M_2\geq \lambda^2h^3|x|^2+\frac{db}{8}\frac{|x|^2}{(t+1)^3}.\\
\end{equation}

\textbf{Estimate of $|\nabla(F-F_0)|$}.

Since
$$F-F_0=2(a^{ij}_\epsilon-a^{ij})(\partial_{ij}\Phi+2\partial_i\Phi\partial_j\Phi),$$
then
\begin{eqnarray*}
\begin{split}
|\nabla(F-F_0)|=&2|(\nabla a^{ij}_\epsilon-\nabla a^{ij})(\partial_{ij}\Phi+
2\partial_i\Phi\partial_j\Phi)\\
&+(a^{ij}_\epsilon-a^{ij})(\nabla\partial_{ij}\Phi+4\partial_i\Phi\nabla\partial_j\Phi)|\\
\leq&2|\nabla a^{ij}_\epsilon-\nabla a^{ij}|(|\partial_{ij}\Phi|+|\partial_i\Phi|^2+|\partial_j\Phi|^2)\\
&+2|a^{ij}_\epsilon-a^{ij}|(|\nabla\partial_{ij}\Phi|+2|\partial_i\Phi|^2+2|\nabla\partial_j\Phi|^2).
\end{split}
\end{eqnarray*}
By (\ref{indexA}) we have
$$|\nabla a^{ij}_\epsilon|\leq\frac{2E}{|x|},~~|a^{ij}_\epsilon-a^{ij}|\leq\frac{E}{|x|},$$
then
\begin{equation}\label{esterror}
\begin{split}
|\nabla(F-F_0)|\leq&\frac{6E}{|x|}(n|\nabla^2\Phi|+2n|\nabla\Phi|^2)+
\frac{2E}{|x|}(n|\nabla^3\Phi|+2n|\nabla\Phi|^2+2n|\nabla^2\Phi|^2)\\
\leq&\frac{n
E}{|x|}(6|\nabla^2\Phi|+2|\nabla^3\Phi|+16|\nabla\Phi|^2+4|\nabla^2\Phi|^2).
\end{split}
\end{equation}
By (\ref{aaa4b}) we have
$$|\nabla(F-F_0)|\leq C(n)E[h^2+\frac{1}{(t+1)^2}]|x|.$$

\textbf{Estimate of $F(x,0)$ and $F(x,1)$.}

By (\ref{choiceofF}) and direct calculations, we have
\begin{eqnarray*}
\begin{split}
F=&2\partial_t\Phi-2a^{ij}\partial_{ij}\Phi-4\langle \mathbf{A}\nabla\Phi,\nabla\Phi\rangle-
\frac{d}{t+1}\\
=&-2\beta\gamma(t+1)^{-\beta-1}|x|^{3/2}-9h^2a^{ij}x_ix_j+
3h[(\frac{1}{2|x|^2}+\frac{4b}{t+1})a^{ij}x_ix_j-a^{ii}]\\
&+\frac{b|x|^2-4b^2a^{ij}x_ix_j}{(t+1)^2}+\frac{\beta}{(t+1)^2}+\frac{2ba^{ii}-d}{t+1},
\end{split}
\end{eqnarray*}
then
\begin{equation*}
\begin{split}
F(x,1)=&-2\beta\gamma2^{-\beta-1}|x|^{3/2}+\frac{b|x|^2-4b^2a^{ij}(x,1)x_ix_j}{4}+
\frac{\beta}{4}+\frac{2ba^{ii}(x,1)-d}{2}\\
\leq&\frac{b|x|^2+\beta}{4}\leq\frac{\beta}{4}(|x|^2+1)\leq\frac{\beta}{2}|x|^2,
\end{split}
\end{equation*}
and
\begin{equation*}
\begin{split}
F(x,0)=&-2\beta\gamma|x|^{3/2}-9\gamma^2(1-2^{-\beta})^2|x|^{-1}a^{ij}(x,0)x_ix_j\\
&+3\gamma(1-2^{-\beta})|x|^{-1/2}[(\frac{1}{2|x|^2}+4b)a^{ij}(x,0)x_ix_j-a^{ii}(x,0)]\\
&+b|x|^2-4b^2a^{ij}(x,0)x_ix_j+\beta+2ba^{ii}-d\\
\geq&-2\beta\gamma|x|^{3/2}-9\gamma^2\Lambda|x|-3\gamma\Lambda|x|^{-1/2}+(b-4b^2\Lambda)|x|^2\\
\geq&-2\beta\gamma|x|^{3/2}(1+\gamma+1)\geq-2\beta|x|^{2}(1+\gamma)^2.
\end{split}
\end{equation*}
Thus we complete the proof of Lemma \ref{pCI1-estimates}.

\subsection{Proof of Proposition \ref{Prop-C2}.}

In this part, we let
$$\Phi=\Psi=\gamma(1-t)R^{2/3}|x|^{4/3}+\psi(t)R^2,$$
and we denote by $c$ absolute constants and $C=C(n,\Lambda,\lambda,M,E)$. We keep in mind that
$$\frac{|x|}{R}\leq1,~~\frac{1}{8}\leq t\leq\frac{7}{8}~~~in~~Q_R.$$

\textbf{Step 1.} Estimate matrix $\mathbf{B}$.

First we estimate the Hessian matrix $D^2\Phi$. Denote
$$g=\gamma(\frac{|x|}{R})^{-2/3}.$$
By direct calculations, we have
$$
D^2\Phi=\frac{4}{3}(1-t)g(I_n-\frac{2x\cdot x^T}{3|x|^2})\geq\frac{4}{9}(1-t)gI_n
\geq cgI_n,
$$
and hence
$$4\mathbf{A}D^2\Phi \mathbf{A}\geq c\lambda^2gI_n.$$

Then we estimate $\partial_l\Phi a^{ki}\partial_ka^{lj}$ and
$\partial_ta^{ij}$.

For any $\xi\in\mathbb{R}^n$,
$$
|\partial_l\Phi a^{ki}\partial_ka^{lj}\xi_i\xi_j|\leq n^2\Lambda\frac{E}{|x|}|\nabla\Phi|\sum_{i,j}|\xi_i||\xi_j|
\leq \frac{n^3\Lambda E}{|x|}|\nabla\Phi||\xi|^2,
$$
Since
$$\nabla\Phi=\frac{4}{3}(1-t)gx,$$
then
\begin{equation}\label{pCI2-graphi}
\frac{1}{6}\gamma R^{2/3}|x|^{1/3}\leq|\nabla\Phi|\leq \frac{4}{3}\gamma R^{2/3}|x|^{1/3},
\end{equation}
and
$$
|\partial_l\Phi a^{ki}\partial_ka^{lj}\xi_i\xi_j|\leq c n^3\Lambda Eg|\xi|^2,
$$
thus
$$\partial_l\Phi a^{ki}\partial_ka^{lj}\geq -cn^3\Lambda EgI_n.$$
Similarly,
$$\partial_ta^{ij}\geq-nMI_n.$$
Consequently,
$$
\mathbf{B}\geq c(\lambda^2-c_1n^3\Lambda E)gI_n-CI_n+H\mathbf{A}.\\
$$
Now we take
$$H=4n^2\varphi(t)Eg,$$
where $\varphi(t)$ is a smooth decreasing function on $[0,1]$ satisfying
$$\varphi(t)=1~~in~[0,\frac{1}{3}],~~~~\varphi(t)=-1~in~[\frac{2}{3},1],$$
$$\varphi(t)>0~~in~[0,\frac{1}{2}),~~~~\varphi(t)<0~~~in~(\frac{1}{2},1].$$
Then
\begin{eqnarray*}
\begin{split}
\mathbf{B}\geq &c(\lambda^2-c_1n^3\Lambda E)g I_n-CI_n-4n^2\Lambda Eg I_n\\
\geq &c(\lambda^2-c_2n^3\Lambda E)g I_n-CI_n,
\end{split}
\end{eqnarray*}
When $E<E_0(n,\Lambda,\lambda)$, and we take $\gamma(n,\Lambda,\lambda,M,E)$ large enough, then
\begin{equation}\label{pCI2-estB}
\mathbf{B}\geq 2c\lambda^2gI_n.
\end{equation}

\textbf{Step 2.} Prove the Carleman inequality.

By (\ref{pCI2-estB}), we have the estimates of the second term of
the left hand side of (\ref{generalCI2}), in fact
\begin{equation}\label{pCI2-3}
\begin{split}
\int_{Q_R}\langle \mathbf{B}\nabla v,\nabla v\rangle e^{2\Phi}dxdt
\geq& 2c\lambda^2\int_{Q_R}ge^{2\Phi}|\nabla v|^2dxdt\\
=&c\lambda^2\int_{Q_R}ge^{2\Phi}|\nabla v|^2dxdt+c\lambda^2\int_{Q_R}g|\nabla w|^2dxdt\\
&+c\lambda^2\int_{Q_R}[g|\nabla\Phi|^2+\nabla g\cdot\nabla\Phi+g\Delta\Phi]w^2dxdt.
\end{split}
\end{equation}
By (\ref{generalCI2}), (\ref{pCI2-3}) and the Cauchy inequality, we have
\begin{equation}\label{pCI2-4}
\begin{split}
&c\lambda^2\int_{Q_R}ge^{2\Phi}|\nabla v|^2dxdt+c\lambda^2\int_{Q_R}g|\nabla w|^2dxdt+\int_{Q_R}M_2w^2dxdt\\
&-\int_{Q_R}w\langle A\nabla(F-F_0),\nabla w\rangle dxdt\leq \int_{Q_R}e^{2\Phi}|Pv|^2dxdt,
\end{split}
\end{equation}
where
$$M_2=c\lambda^2(g|\nabla\Phi|^2+\nabla g\cdot\nabla\Phi+g\Delta\Phi)
+\frac{1}{2}\partial_tF+\frac{1}{2}F(\frac{\partial_tG-\tilde{\Delta}G}{G}-F)+\frac{1}{2}\tilde{\Delta}F_0.$$

We use inequality (\ref{pCI2-4}) to prove Proposition \ref{Prop-C2}.
We also need some estimates which we list in the following lemma.

\begin{Lemma}\label{pCI2-estimates}
There exists a constant $E_0(n,\Lambda,\lambda)$, such that when $E<E_0$, for any \\
$\gamma\geq\gamma_0(n,\Lambda,\lambda,M,E)$, we have
\begin{equation}\label{pCI2-estM2}
M_2\geq c\lambda^2\gamma^3R^2;
\end{equation}
\begin{equation}\label{pCI2-esterror}
|\nabla(F-F_0)|\leq cnE\gamma^2R^{4/3}|x|^{-1/3}.
\end{equation}
\end{Lemma}

We will prove this lemma later.

Then by (\ref{pCI2-esterror}), we have
\begin{eqnarray*}
\begin{split}
|\int_{Q_R}w\langle A\nabla(F-F_0),\nabla w\rangle dx dt|\leq&\Lambda\int_{Q_R}|\nabla(F-F_0)||w||\nabla w|dxdt\\
\leq&cn\Lambda E\int_{Q_R}\gamma^2R^{4/3}|x|^{-1/3}|w||\nabla w|dx dt.
\end{split}
\end{eqnarray*}
Using the Cauchy inequality,
\begin{equation*}
\begin{split}
&|\int_{Q_R}w\langle A\nabla(F-F_0),\nabla w\rangle dxdt|\\
\leq &cn\Lambda E[\int_{Q_R}\gamma^3R^2w^2dxdt+\int_{Q_R}\gamma(\frac{|x|}{R})^{-2/3}|\nabla w|^2dxdt].
\end{split}
\end{equation*}
When $E<E_0(n,\Lambda,\lambda)$, we have
\begin{equation}\label{pCI2-5}
\begin{split}
&|\int_{Q_R}w\langle A\nabla(F-F_0),\nabla w\rangle dxdt|\\
\leq& \frac{1}{2}c\lambda^2[\int_{Q_R}\gamma^3R^2w^2dxdt+\int_{Q_R}\gamma(\frac{|x|}{R})^{-2/3}|\nabla w|^2dxdt]\\
\leq& \frac{1}{2}\int_{Q_R}M_2w^2dxdt+c\lambda^2\int_{Q_R}g|\nabla w|^2dxdt.
\end{split}
\end{equation}
Because of (\ref{pCI2-4}) and (\ref{pCI2-5}), we have
\begin{eqnarray*}
\begin{split}
\int_{Q_R}e^{2\Phi}|Pv|^2dxdt\geq& c\lambda^2\int_{Q_R}ge^{2\Phi}|\nabla v|^2dxdt+\frac{1}{2}\int_{Q_R}M_2w^2dxdt\\
\geq& c\lambda^2\int_{Q_R}e^{2\Phi}(\gamma^3R^2v^2+\gamma|\nabla v|^2)dxdt.
\end{split}
\end{eqnarray*}
Thus we proved Carleman inequality (\ref{C2}).

It remains to prove Lemma \ref{pCI2-estimates}.\\

\textbf{Step 3.} Prove Lemma \ref{pCI2-estimates}.

\textbf{Estimate of $M_2$.}

We estimate the terms of $M_2$ respectively. The leading term of $M_2$ is $h|\nabla\Phi|^2$ and we need pay
attention to two quantities,
$\partial_{t}^2\Phi$ and $\partial_t\Phi(H-2\partial_ia^{ij}\partial_j\Phi)$.\\

\emph{Estimate of the first three terms.}

By (\ref{pCI2-graphi}), we have
$$g|\nabla\Phi|^2\geq c\gamma^3R^2,$$
$$|\nabla g\cdot\nabla\Phi|\leq |\nabla g||\nabla\Phi|\leq c\gamma^2(\frac{|x|}{R})^{-4/3},$$
$$g\Delta\Phi\geq0,$$
then
\begin{equation}\label{t123}
c\lambda^2(g|\nabla\Phi|^2+\nabla g\cdot\nabla\Phi+g\Delta\Phi)\geq c\lambda^2\gamma^3R^2-c\gamma^2(\frac{|x|}{R})^{-4/3}.
\end{equation}

\emph{Estimate of $\frac{1}{2}\partial_tF$.}

Recall (\ref{choiceofF}), then
$$\frac{1}{2}\partial_tF=\partial_t^2\Phi-\partial_ta^{ij}\partial_{ij}\Phi-a^{ij}\partial_{ijt}\Phi
-2\partial_t\langle A\nabla\Phi,\nabla\Phi\rangle-\frac{1}{2}\partial_t H.$$
We estimate them one by one.
$$\partial_t^2\Phi=\psi''R^2\geq-cR^2;$$
$$-\partial_ta^{ij}\partial_{ij}\Phi\geq-C|\nabla^2\Phi|\geq-Cg;$$
$$-a^{ij}\partial_{ijt}\Phi=\frac{4}{3}g(a^{ii}-\frac{2a^{ij}x_ix_j}{3|x|^2})\geq-Cg;$$
\begin{eqnarray*}
\begin{split}
-2\partial_t\langle A\nabla\Phi,\nabla\Phi\rangle
=&-2\partial_ta^{ij}\partial_i\Phi\partial_j\Phi-4a^{ij}\partial_i\Phi\partial_{jt}\Phi\\
\geq&-C|\nabla\Phi|^2+\frac{64}{9}\gamma^2(1-t)(\frac{|x|}{R})^{-4/3}a^{ij}x_ix_j\\
\geq&-C|\nabla\Phi|^2\geq-C\gamma^2R^{4/3}|x|^{2/3};
\end{split}
\end{eqnarray*}
$$-\frac{1}{2}\partial_t H=-2n^2\varphi'(t)Eg\geq0.$$
Combining them together, we have
\begin{equation}\label{t4}
\begin{split}
\frac{1}{2}\partial_tF\geq&-cR^2-Cg-C\gamma^2R^{4/3}|x|^{2/3}\\
\geq&-cR^2-C\gamma^2R^{4/3}|x|^{2/3}.
\end{split}
\end{equation}

\emph{Estimate of $\frac{1}{2}F(\frac{\partial_tG-\tilde{\Delta}G}{G}-F)$.}

First we have
$$\frac{1}{2}F(\frac{\partial_tG-\tilde{\Delta}G}{G}-F)=(\partial_t\Phi-2\langle
A\nabla\Phi,\nabla\Phi\rangle-
a^{ij}\partial_{ij}\Phi-\frac{1}{2}H)(H-2\partial_ia^{ij}\partial_j\Phi).
$$
Since
$$\partial_t\Phi=-\gamma R^{2/3}|x|^{4/3}+\psi'R^2,$$
then
\begin{eqnarray*}
\begin{split}
&\frac{1}{2}F(\frac{\partial_tG-\tilde{\Delta}G}{G}-F)=\psi'R^2(H-2\partial_ia^{ij}\partial_j\Phi)\\
&~~~~~~~~-[\gamma R^{2/3}|x|^{4/3}+2\langle A\nabla\Phi,\nabla\Phi\rangle+a^{ij}\partial_{ij}\Phi
+\frac{1}{2}H](H-2\partial_ia^{ij}\partial_j\Phi) \\
&~~~~\equiv J_1-J_2.
\end{split}
\end{eqnarray*}
Before we estimate $J_1$ and $J_2$, we estimate $2\partial_ia^{ij}\partial_j\Phi$ first.\\
$$|2\partial_ia^{ij}\partial_j\Phi|\leq\frac{2n^2E}{|x|}|\nabla\Phi|,$$
and by (\ref{pCI2-graphi}), we have
$$|2\partial_ia^{ij}\partial_j\Phi|\leq\frac{8}{3}n^2Eg.$$
For $J_1$, we notice that
$$\psi'(t)=0~~in~~[0,\frac{1}{4}]\cup[\frac{1}{3},\frac{2}{3}]\cup[\frac{3}{4},1],$$
so we just need to consider the case when $t\in[\frac{1}{4},\frac{1}{3}]\cup[\frac{2}{3},\frac{3}{4}].$\\
When $t\in[\frac{1}{4},\frac{1}{3}]$, $\psi'\geq0$, $\varphi(t)=1$, $H=4n^2Eg$, then
$$H-2\partial_ia^{ij}\partial_j\Phi\geq0,$$
and thus $J_1\geq0$.\\
When $t\in[\frac{2}{3},\frac{3}{4}]$, $\psi'\leq0$, $\varphi(t)=-1$, $H=-4n^2Eg$, then
$$H-2\partial_ia^{ij}\partial_j\Phi\leq0,$$
and thus $J_1\geq0$.\\
Above all, we have
$$J_1\geq0.$$

For $J_2$,
\begin{eqnarray*}
\begin{split}
J_2\leq&[\gamma R^{2/3}|x|^{4/3}+2\Lambda|\nabla\Phi|^2+C|\nabla^2\Phi|+cn^2Eg]\cdot cn^2Eg\\
\leq&[\gamma R^{2/3}|x|^{4/3}+c\Lambda\gamma^2R^{4/3}|x|^{2/3}+C g]\cdot cn^2Eg\\
\leq&[c\Lambda\gamma^2R^{4/3}|x|^{2/3}+C\gamma R^{2/3}|x|^{4/3}]\cdot cn^2Eg\\
=&cn^2\Lambda E\gamma^3R^2+C\gamma^2R^{4/3}|x|^{2/3}.
\end{split}
\end{eqnarray*}
Combining $J_1$ and $J_2$ together, we obtain
\begin{equation}\label{t5}
\frac{1}{2}F(\frac{\partial_tG-\tilde{\Delta}G}{G}-F)\geq-cn^2\Lambda E\gamma^3R^2-C\gamma^2R^{4/3}|x|^{2/3}.
\end{equation}

\emph{Estimate of $\frac{1}{2}\tilde{\Delta}F_0$.}

Recall that
\begin{eqnarray*}
\begin{split}
F_0=&2\partial_t\Phi-2a^{ij}_\epsilon\partial_{ij}\Phi-4a^{ij}_\epsilon\partial_i\Phi\partial_j\Phi-H\\\
=&-2\gamma R^{2/3}|x|^{4/3}+2\psi'R^2-4n^2\varphi(t)E\gamma R^{2/3}|x|^{-2/3}-
2a^{ij}_\epsilon(\partial_{ij}\Phi+2\partial_i\Phi\partial_j\Phi),
\end{split}
\end{eqnarray*}
then
\begin{equation}\label{t6-1}
\begin{split}
\frac{1}{2}\tilde{\Delta}F_0=&-\gamma R^{2/3}\tilde{\Delta}(|x|^{4/3})-
2n^2\varphi(t)E\gamma R^{2/3}\tilde{\Delta}(|x|^{-2/3})-
\tilde{\Delta}[a^{ij}_\epsilon(\partial_{ij}\Phi+2\partial_i\Phi\partial_j\Phi)]\\
\geq&-C\gamma R^{2/3}|x|^{-2/3}-C\gamma R^{2/3}|x|^{-8/3}-
\tilde{\Delta}[a^{ij}_\epsilon(\partial_{ij}\Phi+2\partial_i\Phi\partial_j\Phi)]\\
\geq&-Cg-\tilde{\Delta}[a^{ij}_\epsilon(\partial_{ij}\Phi+2\partial_i\Phi\partial_j\Phi)],
\end{split}
\end{equation}
and thus it remains to estimate $\tilde{\Delta}[a^{ij}_\epsilon(\partial_{ij}\Phi+
2\partial_i\Phi\partial_j\Phi)]$.\\
By (\ref{indexA}) we have that $|a^{ij}_\epsilon|$, $|\nabla a^{ij}_\epsilon|$
and $|\nabla^2a^{ij}_\epsilon|$ are all bounded,
and it is easy to verify that
\begin{equation}\label{estbb}
|\nabla^k\Phi|\leq C\gamma R^{2/3}|x|^{4/3-k},~~k=1,2,3,4.
\end{equation}
Similarly to (\ref{estmollifier}),
\begin{equation*}
\begin{split}
&|\tilde{\Delta}[a^{ij}_\epsilon(\partial_{ij}\Phi+2\partial_i\Phi\partial_j\Phi)]|\\
\leq& C(|\nabla^2\Phi|+|\nabla^3\Phi|+|\nabla^4\Phi|+|\nabla\Phi|^2+|\nabla^2\Phi|^2+|\nabla^3\Phi|^2),
\end{split}
\end{equation*}
then by (\ref{estbb}),
\begin{equation}\label{t6-2}
|\tilde{\Delta}[a^{ij}_\epsilon(\partial_{ij}\Phi+2\partial_i\Phi\partial_j\Phi)]|
\leq C|\nabla\Phi|^2\leq C\gamma^2R^{4/3}|x|^{2/3}.
\end{equation}
We combine (\ref{t6-1}) and (\ref{t6-2}) and obtain that
\begin{equation}\label{t6}
\frac{1}{2}\tilde{\Delta}F_0\geq-Cg-C\gamma^2R^{4/3}|x|^{2/3}\geq-C\gamma^2R^{4/3}|x|^{2/3}.
\end{equation}
At last, combining (\ref{t123}), (\ref{t4}), (\ref{t5}) and (\ref{t6}), we
have
\begin{eqnarray*}
\begin{split}
M_2\geq &(c\lambda^2-c_3n^2\Lambda E)\gamma^3R^2-C\gamma^2R^{4/3}|x|^{2/3}-cR^2\\
\geq &(c\lambda^2-c_3n^2\Lambda E)\gamma^3R^2-C\gamma^2R^2,
\end{split}
\end{eqnarray*}
When $E<E_0(n,\Lambda,\lambda)$, we have
$$M_2\geq(c\lambda^2\gamma^3-C\gamma^2)R^2\geq c\lambda^2\gamma^3R^2.$$
if $\gamma\geq\gamma_0(n,\Lambda,\lambda,M,E)$ large enough.\\

\textbf{Estimate of $|\nabla(F-F_0)|$}.

Similarly to (\ref{esterror}),
$$|\nabla(F-F_0)|\leq\frac{n E}{|x|}(6|\nabla^2\Phi|+2|\nabla^3\Phi|+16|\nabla\Phi|^2+4|\nabla^2\Phi|^2),$$
then by (\ref{estbb}) we have
$$|\nabla(F-F_0)|\leq\frac{cnE}{|x|}|\nabla\Phi|^2\leq cnE\gamma^2R^{4/3}|x|^{-1/3}.$$
Thus we complete the proof of Lemma \ref{pCI2-estimates}.

\section{Appendix}

\emph{The properties of $\{a^{ij}_\epsilon\}$.}\\

$a^{ij}_\epsilon(x,t)=\int_{\mathbb{R}^n}a^{ij}(x-y,t)\phi_\epsilon(y)dy$,
where $\phi$ is a mollifier and $\epsilon=\frac{1}{2}$, then $\{a^{ij}_\epsilon\}$ satisfy:
\begin{eqnarray*}
\begin{split}
&i)~\lambda|\xi|^2\leq a^{ij}_\epsilon(x,t)\xi_i\xi_j\leq \Lambda|\xi|^2, ~\forall \xi\in\mathbb{R}^n;\\
&ii)~|\nabla a^{ij}_\epsilon(x,t)|\leq M;~~|\nabla a^{ij}_\epsilon(x,t)|\leq\frac{2E}{|x|}~when~|x|\geq1;\\
&iii)~|a^{ij}_\epsilon(x,t)-a^{ij}(x,t)|\leq 2\Lambda;
~~|a^{ij}_\epsilon(x,t)-a^{ij}(x,t)|\leq\frac{E}{|x|}~when~|x|\geq1;\\
&iv)~|\partial_{kl}a^{ij}_\epsilon(x,t)|\leq c(n)M;~~|\partial_{kl}a^{ij}_\epsilon(x,t)|\leq \frac{c(n)E}{|x|}~when~|x|\geq1.
\end{split}
\end{eqnarray*}
\emph{Proof}.\\

$i)$ It is obvious.

$ii)$
$$|\nabla a^{ij}_\epsilon(x,t)|\leq\int_{\mathbb{R}^n}|\nabla a^{ij}(x-y,t)|\phi_\epsilon(y)dy
\leq M\int_{\mathbb{R}^n}\phi_\epsilon(y)dy=M,$$ and when
$|x|\geq1$,
\begin{equation*}
\begin{split}
|\nabla a^{ij}_\epsilon(x,t)|\leq&\int_{\mathbb{R}^n}|\nabla a^{ij}(x-y,t)|\phi_\epsilon(y)dy\\
\leq &\int_{\mathbb{R}^n}\frac{E}{|x-y|}\phi_\epsilon(y)dy\\
\leq &\int_{\mathbb{R}^n}\frac{E}{|x|-\frac{1}{2}}\phi_\epsilon(y)dy
\leq\frac{2E}{|x|}.
\end{split}
\end{equation*}

$iii)$ The first part is obvious. We only need to prove the second
one.
\begin{eqnarray*}
\begin{split}
|a^{ij}_\epsilon(x,t)-a^{ij}(x,t)|&\leq\int_{\mathbb{R}^n}|a^{ij}(x-y,t)-a^{ij}(x,t)|\phi_\epsilon(y)dy\\
&\leq\int_{\mathbb{R}^n}|\nabla a^{ij}(x-\theta y,t)||y|\phi_\epsilon(y)dy, ~~~~~~~(0<\theta<1)\\
\end{split}
\end{eqnarray*}
and when $|x|\geq1$,
$$
|a^{ij}_\epsilon(x,t)-a^{ij}(x,t)|\leq\int_{\mathbb{R}^n}\frac{E}{2|x-\theta y|}\phi_\epsilon(y)dy
\leq\int_{\mathbb{R}^n}\frac{E}{2(|x|-\frac{1}{2})}\phi_\epsilon(y)dy
\leq\frac{E}{|x|}.
$$

$iv)$
\begin{eqnarray*}
\begin{split}
|\partial_{kl}a^{ij}_\epsilon(x,t)|&\leq\int_{\mathbb{R}^n}|
\partial_ka^{ij}(x-y,t)||\partial_l\phi_\epsilon(y)|dy\\
&\leq\epsilon^{-n-1}\int_{\mathbb{R}^n}|\partial_ka^{ij}(x-y,t)||(\partial_l\phi)(\frac{y}{\epsilon})|dy\\
&\leq\frac{M}{\epsilon}||\partial_l\phi||_{L^1}\leq 2M||\nabla\phi||_{L^1},
\end{split}
\end{eqnarray*}
and when $|x|\geq1$,
\begin{eqnarray*}
\begin{split}
|\partial_{kl}a^{ij}_\epsilon(x,t)|&\leq\epsilon^{-n-1}\int_{\mathbb{R}^n}|
\partial_ka^{ij}(x-y,t)||(\partial_l\phi)(\frac{y}{\epsilon})|dy\\
&\leq\epsilon^{-n-1}\int_{\mathbb{R}^n}\frac{E}{|x-y|}|(\partial_l\phi)(\frac{y}{\epsilon})|dy\\
&\leq \frac{2E}{\epsilon|x|}||\partial_l\phi||_{L^1}\leq \frac{4E||\nabla\phi||_{L^1}}{|x|}.
\end{split}
\end{eqnarray*}

Then we finished the proof.

\bigskip

%
%\noindent {\bf Acknowledgments.}
%Authors thank helpful discussions with Professors.
%

\vspace{1cm}

{\small}

\noindent

{\bf Jie Wu}\\
Center for Applied Mathematics, Tianjin University, 92 Weijin Road, Nankai District, Tianjin 300072, P. R. China. \\
Email address: {\bf jackwu@amss.ac.cn}\\

{\bf Liqun Zhang}\\
Hua Loo-Keng Key Laboratory of Mathematics, Chinese Academy of Sciences, Beijing 100190, P. R. China. \\
Email address: {\bf lqzhang@math.ac.cn}

\end{document}